\UseRawInputEncoding
\documentclass[11pt]{amsart}
\usepackage{pieri_review_preamble}
\usepackage{amsthm}
\providecommand{\LL}{\mathcal L}

\title{Two-row Schur times Schubert via parallel classes
of positroid varieties}

\author{Tianyi Yu}
\address{Tianyi Yu, Laboratoire d'Alg\`ebre, de Combinatoire et d'Informatique Math\'ematique, Universit\'e du Qu\'ebec \`a Montr\'eal, Montreal QC, Canada}
\email{yu.tianyi@uqam.ca}

\begin{document}
\begin{abstract}
Schubert structure constants are intersection numbers of Schubert varieties. Finding positive combinatorial formulas for these constants remains a major open problem. We resolve this problem when one of the indexing permutations corresponds to a two-row partition. Our approach uses increasing chains, introduced by Sottile for the one-row case, and builds on the skew Schubert functions studied by Bergeron and Sottile. The key step is to connect
increasing chains and these functions 
with the parallel classes of positroid varieties.
\end{abstract}
\maketitle

\section{Introduction}\label{sec:introduction}
Fix positive integers $k<n$ throughout the paper.
The \definition{complete flag variety} $\Fl(n)$ parametrizes flags
\[
 V_1\subseteq V_2\subseteq\cdots\subseteq V_n=\mathbb C^n,
 \qquad \dim V_i=i.
\]

Its Schubert classes, indexed by permutations in $S_n$, form an integral
basis of $H^*(\Fl(n),\mathbb Z)$. Lascoux and Sch\"utzenberger
introduced the \definition{Schubert polynomials} $\fS_u$ as polynomial
representatives of these classes~\cite{LS:Schubert}.
Write $c^w_{u,v}$ for the coefficient of $\fS_w$ in $\fS_u\fS_v$.
For $u,v,w\in S_n$, these \definition{Schubert structure constants}
are intersection numbers and hence nonnegative integers.
A central problem in Schubert calculus is to give them positive combinatorial formulas. Several rules
are known in special cases~\cite{Kogan:RC,Lenart:growth,Huang:separated,KZ:separated}.

A classical special case comes from the Grassmannian $\Gr(k,n)$,
which parametrizes $k$-dimensional subspaces of $\mathbb C^n$.
Its Schubert classes are indexed by partitions
$\lambda=(\lambda_1,\ldots,\lambda_k)$ with
$n-k\ge\lambda_1\ge\cdots\ge\lambda_k\ge0$.
These partitions correspond to \definition{$k$-Grassmannian permutations},
namely $w\in S_n$ such that $w(i)<w(i+1)$ for every
$i\ne k$.
The correspondence sends $w$ to $(w(k)-k,\ldots,w(1)-1)$.
We identify a $k$-Grassmannian permutation with its partition and omit
trailing zeros. Its Schubert polynomial is the \definition{Schur polynomial}
$\fS_\lambda=s_\lambda(x_1,\ldots,x_k)$.
When $u,v,w$ are all $k$-Grassmannian, the constants $c^w_{u,v}$ are
\definition{Littlewood--Richardson coefficients}, which admit
several positive combinatorial formulas.

A natural extension is to study $c^w_{u,\lambda}$ for arbitrary $u,w\in S_n$
and a partition $\lambda$.
Three fundamental results illuminate this problem.
\begin{itemize}
\item Monk's formula~\cite{Monk} states that $c^w_{u,(1)}=1$ when
$u$ is covered by $w$ in the $k$-Bruhat order (Definition~\ref{D: k-Bruhat}), and is zero otherwise.
\item 
Sottile's Pieri rule~\cite{Sot} treats the one-row case $\lambda = (\lambda_1)$ using
\definition{increasing chains} (Definition~\ref{D: inc chain}).
Write $u\xrightarrow{k}w$ when such a chain exists. 
Then Sottile showed
$c^w_{u,(\lambda_1)}=1$ if $u\xrightarrow{k}w$ and
$\lambda_1=\ell(w)-\ell(u)$, and is $0$ otherwise.
Sottile also gave a rule for the hook case
$\lambda=(\lambda_1,1,\ldots,1)$.

\item M\'esz\'aros, Panova, and Postnikov~\cite{MPP:SchurSchubert} 
used the quadratic algebra of Fomin and
Kirillov~\cite{FK:quadratic} and Postnikov's algebraic Pieri
formulas~\cite{Postnikov:quantumPieri}.
They proved the existence of a nonnegative expansion in that algebra when $\lambda$ has the form $(\lambda_1,2,1,\ldots,1)$.
\end{itemize}

M\'esz\'aros, Panova, and Postnikov~\cite[\S6]{MPP:SchurSchubert}
proposed the two-row case $\lambda=(\lambda_1,\lambda_2)$
as a next step.
Recently, Tao and Zheng~\cite{TaoZheng:twoRow}
gave bounds on these coefficients. 
We give an explicit positive combinatorial formula 
for these two-row coefficients.
Let $\LL_{u,w}=\{u(j):j>k,\ u(j)=w(j)\}$.
Define a set partition $\PP_{u,w}$ of $[n]\setminus\LL_{u,w}$:
Numbers $a<b$ in $[n]\setminus\LL_{u,w}$ lie in the same class exactly when
$wu^{-1}((a,b])=[a,b)$.
We express
$c^w_{u,(\lambda_1,\lambda_2)}$ in terms of Kostka numbers.

\begin{thm}\label{T: two row coefficients}
For $u\le_k w$, we have
\begin{equation}\label{EQ: introduction numerical bound}
 \ell(w)-\ell(u)\ge n-2k-|\LL_{u,w}|+|\PP_{u,w}|.
\end{equation}
If the inequality is strict, every two-row coefficient
$c^w_{u,(\lambda_1,\lambda_2)}$ is zero. Otherwise,
\begin{equation}\label{EQ: two row coeff formula}
 c^w_{u,(\lambda_1,\lambda_2)}
 =K_{(\lambda_1-d,\lambda_2-d),\alpha},
 \qquad \text{where }
 d=|\PP_{u,w}|-k,\quad
 \alpha=(|C|-1:C\in\PP_{u,w}).
\end{equation}
\end{thm}
Here $\lambda_1\ge\lambda_2\ge0$, and $K_{\mu,\alpha}$ denotes the
\definition{Kostka number}, which counts semistandard Young tableaux
of shape $\mu$ and weight $\alpha$. By convention, it is zero if
$\mu$ has a negative part or $|\mu|\ne|\alpha|$.

\subsection*{Consecutive increasing chains}
Our strategy is to study pairs of increasing chains
$u\xrightarrow{k}v\xrightarrow{k}w$.
One motivation comes from the Jacobi--Trudi identity
\[
 s_{(a,b)}=h_a h_b-h_{a+1}h_{b-1},
\]
which, together with the Pieri rule, expresses $c^w_{u,(a,b)}$ as a
difference of two counts of such pairs. Here $h_j=s_{(j)}$ is the
complete homogeneous symmetric polynomial, with $h_0=1$ and $h_{-1}=0$.
This Jacobi--Trudi approach is used by Tao and
Zheng~\cite{TaoZheng:twoRow} and is proposed for
two-row Dunkl evaluations by M\'esz\'aros, Panova, and
Postnikov~\cite{MPP:SchurSchubert}.

The role of chains in Schubert multiplication was developed further by
Bergeron and Sottile~\cite{BS:BruhatGeometry,BS:monoid}.
Related chain constructions occur in formulas for skew Schubert
polynomials~\cite{LenartSottile:skew} and in the $K$-theoretic
Pieri rule~\cite{LenartSottile:KPieri}.
More generally, for $m\ge1$, let
\begin{equation}\label{EQ: introduction chain family}
 \C(u,w,m)=\{(v_0,\ldots,v_m):
 u=v_0\xrightarrow{k}v_1\xrightarrow{k}\cdots
 \xrightarrow{k}v_m=w\}.
\end{equation}
Their generating function is the
\definition{skew Schubert function} of Bergeron and
Sottile~\cite{BS:skewSchubert}:
\begin{equation}\label{EQ: ordinary chain generating function}
 F_{u,w}(x_1,\ldots,x_m)
 =\sum_{(v_0,\ldots,v_m)\in\C(u,w,m)}
   \prod_{i=1}^m x_i^{\ell(v_i)-\ell(v_{i-1})}.
\end{equation}
Bergeron and Sottile showed that
these polynomials are symmetric and satisfy
\begin{equation}\label{EQ: skew Schubert expansion}
 F_{u,w}(x_1,\ldots,x_m)
 =\sum_\lambda c^w_{u,\lambda}s_\lambda(x_1,\ldots,x_m).
\end{equation}
Thus $F_{u,w}(x_1,x_2)$ records exactly the two-row
coefficients $c^w_{u,(\lambda_1,\lambda_2)}$.

\begin{figure}[ht]
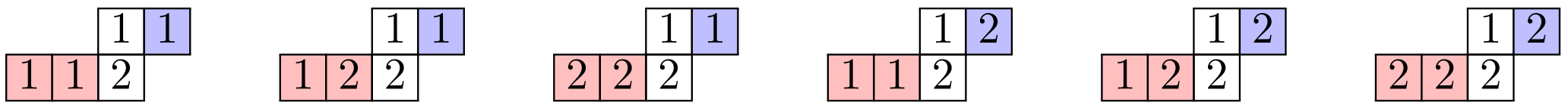

\centering
\[
\ytableausetup{boxsize=1em}
\begin{ytableau}
\none & \none & 1 & *(blue!25)1 \\
*(red!25) 1 & *(red!25)1 & 2
\end{ytableau}
\qquad
\begin{ytableau}
\none & \none & 1 & *(blue!25)1 \\
*(red!25)1 & *(red!25)2 & 2
\end{ytableau}
\qquad
\begin{ytableau}
\none & \none & 1 & *(blue!25)1 \\
*(red!25)2 & *(red!25)2 & 2
\end{ytableau}
\qquad
\begin{ytableau}
\none & \none & 1 & *(blue!25)2 \\
*(red!25)1 & *(red!25)1 & 2
\end{ytableau}
\qquad
\begin{ytableau}
\none & \none & 1 & *(blue!25)2 \\
*(red!25)1 & *(red!25)2 & 2
\end{ytableau}
\qquad
\begin{ytableau}
\none & \none & 1 & *(blue!25)2 \\
*(red!25)2 & *(red!25)2 & 2
\end{ytableau}
\]
\caption{The six semistandard tableaux of shape $(4,3)/(2)$
with entries in $\{1,2\}$. Their generating function is
$s_{(4,3)/(2)}(x_1,x_2)=x_1x_2h_2(x_1,x_2)h_1(x_1,x_2)$.
These six tableaux correspond to $[0,2]\times[0,1]$:
we can independently decide whether we want zero, one,
or two $1$s in the red region, and zero or one $1$
in the blue region.}
\label{F: SSYT}
\end{figure}

Skew Schubert functions recover a well-studied
family of symmetric polynomials.
When $u,w$ are both $k$-Grassmannian and correspond to
$\mu\subseteq\lambda$, we have $F_{u,w}=s_{\lambda/\mu}$, the
\definition{skew Schur polynomial}.
If $s_{\lambda/\mu}(x_1,x_2)$ is nonzero, it has the form
$(x_1x_2)^d h_\alpha(x_1,x_2)$, where
$h_\alpha:=\prod_i h_{\alpha_i}$.
Figure~\ref{F: SSYT} illustrates this phenomenon.

Our main result gives an analogous description for arbitrary
$F_{u,w}(x_1,x_2)$. For $(u,v,w)\in\C(u,w,2)$, define
\begin{equation}\label{EQ: introduction I}
 I(u,v,w)=\{r\in[n]:vu^{-1}(r)=wu^{-1}(r)<r\}.
\end{equation}
When $\C(u,w,2)$ is nonempty, we will show that the image of $I$
is precisely the set of ideals of a poset constructed from $\PP_{u,w}$.
\begin{defn}
The poset $\Gamma_{u,w}$ is the disjoint union of the chain posets
\[
 \Gamma_C:\quad uw^{-1}(c_0)\prec uw^{-1}(c_1)\prec\cdots
 \prec uw^{-1}(c_{L-1}),
 \qquad C=\{c_0<\cdots<c_L\}\in\PP_{u,w}.
\]
A \definition{(lower) ideal} is a subset closed under taking smaller
elements in this poset.
\end{defn}

\begin{theorem}\label{T: main}
For $u\le_k w$, equality in~\eqref{EQ: introduction numerical bound}
holds if and only if $\C(u,w,2)\ne\varnothing$.
If so, the map $I$ is a bijection from $\C(u,w,2)$ to the
ideals of $\Gamma_{u,w}$, with inverse $\Chain$ defined in
Section~\ref{S: reconstruction}.

For $(u,v,w)\in\C(u,w,2)$, let
$d,\alpha$ be as in~\eqref{EQ: two row coeff formula} and put
$S=I(u,v,w)$. Then
\[
 \ell(v)-\ell(u)=d+|S|,\qquad
 \ell(w)-\ell(v)=d+|\alpha|-|S|.
\]
Thus, if equality in~\eqref{EQ: introduction numerical bound}
holds,
\begin{equation}\label{EQ: introduction product}
 F_{u,w}(x_1,x_2)=(x_1x_2)^d h_\alpha(x_1,x_2).
\end{equation}
\end{theorem}
Expanding this product in Schur polynomials proves the coefficient
formula in Theorem~\ref{T: two row coefficients}.

\begin{exa}\label{Ex: three parallel pairs}
Consider $n = 6$ and $k = 2$.
Take $u=[3,1,6,5,4,2]$ and $w=[6,4,5,2,3,1]$ 
from $S_6$.
Then $\rho=wu^{-1}=[4,1,6,3,2,5]$ and
\[
 \LL_{u,w}=\varnothing,\qquad
 \PP_{u,w}=\bigl\{\{1,2\},\{3,4\},\{5,6\}\bigr\}.
\]
For instance, $1$ and $2$ are in the same class
since $\rho(\{2\})=\{1\}$.
Here $\ell(u)=8$, $\ell(w)=13$, $d=3 - 2 = 1$, and $\alpha=(1,1,1)$.
Inequality~\eqref{EQ: introduction numerical bound} is an equality, so
\[
 F_{u,w}(x_1,x_2)
 =x_1x_2(x_1+x_2)^3
 =s_{(4,1)}(x_1,x_2)+2s_{(3,2)}(x_1,x_2).
\]
In particular, $c^w_{u,(3,2)}=2$ and $c^w_{u,(4,1)}=1$.
\end{exa}

\subsection*{Positroid varieties}

The set $\LL_{u,w}$ and set partition $\PP_{u,w}$ have natural geometric
interpretations in terms of positroid varieties.
Suppose $u \leq_k w$.
Let $X_w$ be the Schubert variety indexed by $w$ and $X^u$ the opposite
Schubert variety indexed by $u$. 
Their intersection is a Richardson
variety. Its projection from $\Fl(n)$ to $\Gr(k,n)$ is the
\definition{positroid variety} $\Pi_{u,w}$, of dimension
$\ell(w)-\ell(u)$~\cite[Theorem~5.9]{KLS:positroids}.
The coefficients $c^w_{u,\lambda}$ can be recovered as intersection
numbers of $\Pi_{u,w}$ with Schubert varieties in $\Gr(k,n)$~\cite[Proposition~7.10]{Pawlowski:positroids}.
Their combinatorics and geometry have also been studied through total
positivity~\cite{Postnikov:total}, noncrossing
partitions~\cite{ARW:noncrossing}, and cluster
algebras~\cite{GalashinLam:cluster}.

A general point of $\Pi_{u,w}$ can be viewed
as the row span of a $k\times n$ matrix.
An integer $c$ is a \definition{loop} when column $c$ is zero.
The set of loops is exactly $\LL_{u,w}$.
Distinct nonloops $a,b$ are \definition{parallel}, written $a\parallel b$,
when their columns are proportional. This relation partitions
$[n]\setminus\LL_{u,w}$ into \definition{parallel classes}.
This set partition is exactly $\PP_{u,w}$.
Neither $\LL_{u,w}$ nor $\PP_{u,w}$ depends on the chosen general
point or its representing matrix.
Theorem~\ref{T: two row coefficients} thus gives a vanishing test for
$c^w_{u,(\lambda_1,\lambda_2)}$ using the dimension, the loops,
and the number of parallel classes.
When equality holds in~\eqref{EQ: introduction numerical bound},
the sizes of the parallel classes determine every two-row coefficient.

\subsection*{Double weights}
There is also a refinement of $F_{u,w}(x_1,x_2)$
involving $y$-variables $\y=(y_1,\ldots,y_n)$.
Following Chou and Yu~\cite{CY:doubleEG}, define
\[
 F'_{u,w}(x_1,x_2;\y)
 =\sum_{(u,v,w)\in\C(u,w,2)}
 \prod_{a\in\LL_{u,v}}(x_1-y_a)
 \prod_{b\in\LL_{v,w}}(x_2-y_b).
\]

Chou and Yu~\cite[Proposition~2.29]{CY:doubleEG}, using Samuel's
Pieri identity~\cite[Theorem~7.1]{Samuel:Pieri}, showed that this
polynomial is symmetric in $x_1,x_2$. 
By Lemma~\ref{L: double weight specialization}, it is related to
$F_{u,w}(x_1,x_2)$ by
\[
 F'_{u,w}(x_1,x_2;\mathbf0)
 =(x_1x_2)^{n-k}
 F_{u,w}(x_1^{-1},x_2^{-1}).
\]
We can also describe this polynomial explicitly.

\begin{theorem}\label{T: double two row}
If equality holds in~\eqref{EQ: introduction numerical bound}, then
\begin{equation}\label{EQ: introduction double product}
\begin{aligned}
 F'_{u,w}(x_1,x_2;\y)
 ={}\prod_{a\in\LL_{u,w}}(x_1-y_a)(x_2-y_a)
 \quad\cdot\prod_{C\in\PP_{u,w}}
 \frac{\prod_{c\in C}(x_1-y_c)-\prod_{c\in C}(x_2-y_c)}{x_1-x_2}.
\end{aligned}
\end{equation}

If the inequality is strict, then
$F'_{u,w}=0$.
\end{theorem}

\subsection*{Organization of the proof}
Sections~\ref{S:increasing-chains} and~\ref{S:positroids} give the chain
and positroid preliminaries.
Section~\ref{S: reconstruction} gives an algorithm that reconstructs
$(u,v,w)$ from $I(u,v,w)$ and proves the injectivity of $I$.
Section~\ref{S: inequalities} establishes the inequality~\eqref{EQ: introduction numerical bound} and constructs two
increasing chains when equality holds.
Section~\ref{S: parallel classes} proves the converse: when
$\C(u,w,2)\ne\varnothing$, equality holds
in~\eqref{EQ: introduction numerical bound}. 
Section~\ref{S: image of I} shows
that each $I(u,v,w)$ is an ideal of $\Gamma_{u,w}$,
while Section~\ref{S: surjectivity} shows that
every ideal of $\Gamma_{u,w}$ lies in the image.
Finally, Section~\ref{S: main proofs} assembles these results
to prove the theorems stated in the introduction.

\section*{Acknowledgments}
The author thanks Nantel Bergeron, Cristian Lenart, Frank Sottile, and Alex Yong
for helpful conversations.

\section*{Use of artificial intelligence}
The author used ChatGPT (OpenAI) to assist with language editing,
refining proofs, generating examples, and analyzing computational data.
In particular, ChatGPT helped identify the poset $\Gamma_{u,w}$
by recognizing patterns in computational data.
The author takes full responsibility for the mathematical content
and any remaining errors.

\section{Background and preliminaries on increasing chains}\label{S:increasing-chains}
We review $k$-Bruhat order and the Pieri rule for increasing chains.
We then describe these chains by decreasing cycles and ordered-path graphs.

We regard permutations in $S_n$ as functions on $[n]:=\{1,\dots,n\}$ and represent $u\in S_n$ by its \definition{one-line notation} $[u(1),\dots,u(n)]$.
For $a<b$, write $t_{a,b}\in S_n$ for the transposition interchanging $a$ and $b$.
Left multiplication swaps values, whereas right multiplication swaps positions:
$t_{a,b}u$ is obtained by interchanging the values $a,b$ in $u$, while $ut_{a,b}$ is obtained by interchanging its entries in positions $a,b$.
For a set $S \subseteq [n]$, we let $u(S)$ be the image of $S$ under $u$.
The \definition{length} of $u\in S_n$ is
\begin{equation}\label{EQ: length}
 \ell(u)=|\{(i,j):1\le i<j\le n,\ u(i)>u(j)\}|.
\end{equation}

\subsection{\texorpdfstring{$k$}{k}-Bruhat order}
\begin{defn}\label{D: k-Bruhat}
Let $v,w\in S_n$. We say that $w$ covers $v$ in the \definition{Bruhat order}, written $v\lessdot w$, if $w=vt_{i,j}$ for some $i<j$ and $\ell(w)=\ell(v)+1$.
It is a cover in the \definition{$k$-Bruhat order}, written $v\lessdot_k w$, if additionally $i\le k<j$.
The reflexive and transitive closure of these covers is denoted by $\le_k$.
Finally, we say that $t_{a,b}$ is a \definition{$k$-cover} on $u$ if $u\lessdot_k t_{a,b}u$.
\end{defn}

We slightly generalize the notion of $k$-cover.
A transposition $t_{a,b}$ is a \definition{$k$-swap} on $u$ if $a<b$
and $u^{-1}(a) \leq k < u^{-1}(b)$.
A $k$-swap acting on the left of $u$ increases the length, but need not increase it by one.
It increases the length by $1$ if and only if
it is a $k$-cover on $u$.

We record two useful tests. The first tests whether a transposition gives a cover, while the second tests whether $u\le_k w$.
\begin{lem}\label{L:cover-test}
Let $t_{a,b}$ be a $k$-swap on $v$, and put
$i=v^{-1}(a)$ and $j=v^{-1}(b)$. Then
\[
 \ell(t_{a,b}v)-\ell(v)
 =1+2|\{s:i<s<j,\ a<v(s)<b\}|.
\]
In particular, the $k$-swap $t_{a,b}$ 
is a $k$-cover on $v$ if and only if there is no such $s$.
\end{lem}
\begin{proof}
This follows from~\eqref{EQ: length}.
\end{proof}

\begin{prop}[{\cite[Proposition 3.1]{BS:monoid}}]\label{P:k-Bruhat-endpoints}
For $u,w\in S_n$, we have $u\le_k w$ if and only if
\[
 u(i)\le w(i)\text{ for }i\le k,\qquad
 u(j)\ge w(j)\text{ for }j>k, \text{ and}
\]
\[
 u(i)<u(j)\ \Longrightarrow\ w(i)<w(j)
 \text{ for }i<j\le k\text{ or }k<i<j.
\]
\end{prop}

\subsection{Increasing chains}
We now introduce the saturated chains that appear in the Pieri rule.
\begin{defn}\label{D: inc chain}
Consider a saturated chain
\[
 v_0\lessdot_k v_1\lessdot_k\cdots\lessdot_k v_m.
\]
Write $v_s=t_{a_s,b_s}v_{s-1}$ with $a_s<b_s$ for $1\le s\le m$.
The chain is an \definition{increasing chain} if $a_1<\cdots<a_m$.
In other words, the smaller value exchanged at each step increases.
It is a \definition{right-distinct chain} if the right positions involved in its transpositions are distinct.
\end{defn}

Every increasing chain is right-distinct.
Moreover, there is an increasing chain from $u$ to $w$ if and only if there is a right-distinct chain from $u$ to $w$.
When these chains exist, the increasing chain is unique, although the 
right-distinct chain might not be unique.
Write $u\xrightarrow{k}w$ in this case.
We consider the empty chain from $u$
to itself as an increasing chain, 
so $u\xrightarrow{k}u$.

\begin{exa}
\label{Ex: inc chain}
Let $n = 7$ and $k = 3$.
Take $u = [3,5,1,6,7,4,2]$ and $w = [4,5,7,2,6,3,1]$.
Then the right-distinct chains from $u$ to $w$ are depicted as follows: 
For readability, brackets and commas are omitted when we write a permutation.
We also put a bar between the first $k$
numbers and the last $n-k$ numbers. 
Every $t_{a,b}$ indicates a left action, so 
it swaps the value $a$ and the value $b$. 
\[
\begin{gathered}
351|6742\xrightarrow{t_{1,2}}352|6741\xrightarrow{t_{2,6}}356|2741\xrightarrow{t_{3,4}}456|2731\xrightarrow{t_{6,7}}457|2631\\[4pt]
351|6742\xrightarrow{t_{1,2}}352|6741\xrightarrow{t_{2,6}}356|2741\xrightarrow{t_{6,7}}357|2641\xrightarrow{t_{3,4}}457|2631\\[4pt]
351|6742\xrightarrow{t_{1,2}}352|6741\xrightarrow{t_{3,4}}452|6731\xrightarrow{t_{2,6}}456|2731\xrightarrow{t_{6,7}}457|2631\\[4pt]
351|6742\xrightarrow{t_{3,4}}451|6732\xrightarrow{t_{1,2}}452|6731\xrightarrow{t_{2,6}}456|2731\xrightarrow{t_{6,7}}457|2631
\end{gathered}
\]
The first is the unique increasing chain from $u$ to $w$.
Since such a chain exists, we write $u \xrightarrow{3}w$.
\end{exa}

The Pieri rule connects increasing chains with the Schubert structure
constants $c^w_{u,v}$. When $v$ is the Grassmannian permutation
corresponding to the one-row partition $(r)$, it gives the following
explicit description of $c^w_{u,(r)}$.
\begin{thm}[{\cite[Theorem 1]{Sot}}]\label{T:Pieri-background}
For $r\ge0$,
\[
 c^w_{u,(r)}=
 \begin{cases}
 1,&u\xrightarrow{k}w\text{ and }\ell(w)-\ell(u)=r,\\
 0,&\text{otherwise}.
 \end{cases}
\]
\end{thm}

For later arguments, we write the chain in
Definition~\ref{D: inc chain} more compactly as
\begin{equation}\label{EQ: trans inc chain}
 v_m=t_{a_m,b_m}\cdots t_{a_2,b_2}t_{a_1,b_1}v_0.
\end{equation}
The rightmost transposition acts first and each $t_{a_i, b_i}$
is a $k$-cover on the permutation $t_{a_{i-1},b_{i-1}}\cdots t_{a_1,b_1}v_0$.

\begin{exa}
\label{Ex: swaps}
Continuing Example~\ref{Ex: inc chain},
we may write the right-distinct chains as 
\[
\begin{aligned}
w&=t_{6,7}t_{3,4}t_{2,6}t_{1,2}u\\
w&=t_{3,4}t_{6,7}t_{2,6}t_{1,2}u\\
w&=t_{6,7}t_{2,6}t_{3,4}t_{1,2}u\\
w&=t_{6,7}t_{2,6}t_{1,2}t_{3,4}u. \qedhere
\end{aligned}
\]
\end{exa}

\subsection{Decreasing cycles}
The transpositions of an increasing chain can be grouped into disjoint
decreasing cycles. We recall this description from~\cite{BS:skewSchubert}.
A \definition{cycle} $(a_0,a_1,\ldots,a_s)$ sends $a_i$ to $a_{i+1}$ for $0\le i<s$ and sends $a_s$ to $a_0$.
Its \definition{absolute length} is $s$.
In particular, $(a_0)$ is a cycle of absolute length zero.
We say the cycle is \definition{decreasing} if $a_0>\cdots>a_s$.
With this convention,
\begin{equation}\label{EQ: cycle to trans}
 (a_0,a_1,\ldots,a_s)
 =t_{a_1, a_0} t_{a_2, a_1} \cdots t_{a_{s}, a_{s-1}}.
\end{equation}

\begin{prop}[{\cite[Proposition 1.5]{BS:skewSchubert}}]
\label{P: inc chain cycles}
Take $u,w\in S_n$. Then $u\xrightarrow{k}w$ if and only if there are $k$ decreasing cycles $c_1,\ldots,c_k$ satisfying the following conditions.
\begin{itemize}
\item The cycles involve distinct numbers.
\item The cycle $c_i$ ends with $u(i)$ for $i \in [k]$.
\item We have $w=c_1\cdots c_ku$.
\item The sum of their absolute lengths is $\ell(w)-\ell(u)$.
\end{itemize}
These cycles, if they exist, are unique. 
Expanding each decreasing cycle as in~\eqref{EQ: cycle to trans} gives a right-distinct chain in the style~\eqref{EQ: trans inc chain}.
\end{prop}

We record a ``non-crossing'' property of these cycles. 
\begin{lem}[{\cite[proof of Corollary 5.3]{BS:skewSchubert}}]
\label{L: cycle noncross}
Suppose the cycles $c_1, \dots, c_k$ in Proposition~\ref{P: inc chain cycles} exist. Then they are non-crossing, meaning that
there do not exist $a < b < c < d$ such that $a$ and $c$ are in one cycle
while $b$ and $d$ are in another cycle. 
\end{lem}

\begin{exa}
For $u,w$ in Example~\ref{Ex: inc chain}, the three decreasing
cycles are
\[
 c_1=(4,3),\qquad c_2=(5),\qquad c_3=(7,6,2,1).
\]
They end with $u(1)=3$, $u(2)=5$, and $u(3)=1$, respectively, and
\[
 w=(4,3)(5)(7,6,2,1)u,
 \qquad 1+0+3=\ell(w)-\ell(u)=4.
\]
Expanding these cycles gives
\[
 w=t_{3,4}t_{6,7}t_{2,6}t_{1,2}u,
\]
which is the second right-distinct chain displayed in Example~\ref{Ex: swaps}.
\end{exa}

\begin{rem}
\label{R: force k-swap to be cover}
If $w=t_{a_s,b_s}\cdots t_{a_1,b_1}u$ and $t_{a_i,b_i}$
is a $k$-swap on $t_{a_{i-1}, b_{i-1}}\cdots t_{a_1, b_1}u$,
then $\ell(w)-\ell(u)\ge s$.
If $\ell(w)-\ell(u)=s$, every swap must be a $k$-cover.
Thus the last condition in Proposition~\ref{P: inc chain cycles} can 
be replaced by the condition that the sum of the absolute lengths is at least $\ell(w)-\ell(u)$.
   
\end{rem}

\subsection{Increasing chain graphs}
A decreasing cycle may instead be recorded by a path through its
numbers in increasing order. This leads to an encoding of increasing chains by graphs. 
For a graph $G$, we denote its vertex set and edge set by 
$V(G)$ and $E(G)$.

We consider only undirected graphs with $V(G) \subseteq [n]$.
For us, an edge is a pair $(a,b)$ where $a,b \in V(G)$ 
and $a < b$.
We say $G$ is an \definition{ordered path}
if $$V(G) = \{a_0 < a_1 < \cdots < a_s\} \text{ and }E(G) = \{(a_0, a_1), (a_1, a_2), \dots, (a_{s-1}, a_s)\}.$$
We define the \definition{head} and \definition{tail} of such a path
as $a_0$ and $a_s$ respectively.
Notice that a graph with one vertex
is considered as an ordered path.
We call such a graph a \definition{singleton}.
The head and tail of a singleton are both that vertex.

We say $G$ \definition{decomposes into ordered paths}
if each connected component of $G$ is an ordered path. 
Let $\Head(G)$ and $\Tail(G)$ be the set of heads and tails respectively
over all paths in $G$.

\begin{rem}
\label{R: Two graphs}
Let $G$ and $G'$ be two graphs that decompose into
ordered paths. 
If $E(G)=E(G')$, then $G$ and $G'$ agree after ignoring their singletons. 
If we also have $\Head(G)=\Head(G')$, then $G=G'$.
\end{rem}

Define the permutation $\pi(G) \in S_n$ as a product of cycles.
Each path in $G$ contributes a decreasing cycle by writing its
vertices decreasingly. 
The effect of $\pi(G)$ has a more explicit description.
\begin{itemize}
\item For each ordered path with head
$h$ and tail $t$, $\pi(G)$ sends $h$ to $t$.
\item For each edge $(a,b)$ in $G$, $\pi(G)$ sends $b$
to $a$.
\item All other numbers are fixed by $\pi(G)$. 
\end{itemize}

\begin{cor}
\label{C: Graph language}
Take $u,w\in S_n$. Then $u\xrightarrow{k}w$ if and only if there 
exists a graph $G$ decomposing into $k$ ordered paths
such that $\Head(G) = u([k])$, $w = \pi(G) u$ and $|E(G)| = \ell(w)-\ell(u)$.
This graph $G$, if it exists, is unique.
In addition, the last condition may be replaced by
$|E(G)| \geq \ell(w)-\ell(u)$.
\end{cor}

In this case, we say $G$ is the \definition{increasing chain graph}
of $u$ and $w$.

\begin{proof}
This is a translation of Proposition~\ref{P: inc chain cycles}:
each cycle becomes an ordered path using the numbers in the cycle.
 
The last statement follows from Remark~\ref{R: force k-swap to be cover}.
\end{proof}

\begin{exa}
The increasing chain graph for Example~\ref{Ex: inc chain} can be drawn as
\[
 3\! -\! 4,\qquad 5,\qquad 1\! -\! 2\! -\! 6\! -\! 7.
\]
Explicitly, $V(G)=[7]$ and $E(G) = \{(1,2),(2,6),(3,4),(6,7)\}$.
We have
\[
\begin{split}
 \Head(G)=\{1,3,5\}=u([3]) \text{ and }\Tail(G)=\{4,5,7\}=w([3]).
\end{split}
\]
Its associated permutation is
$\pi(G)=(4,3)(5)(7,6,2,1)$.
\end{exa}

We record some facts about the increasing chain graph.
\begin{lem}\label{L: first graph position values}
Suppose $G$ is the increasing chain graph from $u$ to $w$.
The graph $G$ has the following vertices and edges:
\[
 V(G) = u([k]) \sqcup \{u(j): u(j) > w(j)\} \text{ and } E(G)=\{(w(j),u(j)): u(j)>w(j)\}.
\]
Also, $\Tail(G) = w([k])$.
\end{lem}
\begin{proof}
For each edge $(a,b) \in E(G)$,
$\pi(G)$ sends $b$ to $a$.
The description of $E(G)$ then follows from
$w = \pi(G) u$.
By definition, $u([k])  = \Head(G) \subseteq V(G)$.
By the description of $E(G)$, 
$\{u(j): u(j) > w(j)\} \subseteq V(G)$.
Now take $j > k$ with $u(j) = w(j)$.
If $u(j) \in V(G)$, the identity of $E(G)$
forces this vertex to be a singleton, 
contradicting $\Head(G) = u([k])$.
Thus, such $u(j) \notin V(G)$,
proving the statement about $V(G)$.

Finally, take $t \in \Tail(G)$.
The identity of $E(G)$ says $t$ 
is not $w(j)$ for any $j$ with $w(j) < u(j)$.
In addition, $t \in V(G)$ implies
that $t$ is not $w(j)$ with $j > k$ and $w(j) = u(j)$.
Thus, $t \in w([k])$.
This shows $\Tail(G) \subseteq w([k])$.
Since $\Head(G) = u([k])$, 
$\Tail(G)$ also has size $k$, so 
$\Tail(G)=  w([k])$.
\end{proof}

\begin{cor}\label{C: pair swap positions}
Suppose $u\xrightarrow{k}w$, with increasing chain graph $G$.
For every $i>k$, we have $w(i)\notin\Tail(G)$, and exactly one
of the following holds:
\begin{enumerate}
\item $u(i)=w(i)$ and $w(i)\notin V(G)$;
\item $w(i) < u(i)$ and $(w(i),u(i))\in E(G)$.
\end{enumerate}
\end{cor}
\begin{proof}
This follows from Lemma~\ref{L: first graph position values}.
\end{proof}

\begin{cor}
\label{C: Edge and increasing swaps}
Suppose $u\xrightarrow{k}w$, and let
$w=t_{a_m,b_m}\cdots t_{a_2,b_2}t_{a_1,b_1}u$
be a right-distinct chain, with $a_t<b_t$ for every $t$.
Then its increasing chain graph $A$
has edges $E(A) = \{(a_i, b_i) : i \in [m]\}$.
\end{cor}
\begin{proof}
Each position greater than $k$ is used at most once
in the chain.
If the transposition
$t_{a_t,b_t}$ uses position $q>k$, then $u(q)=b_t$ and $w(q)=a_t$.
The assertion follows from Lemma~\ref{L: first graph position values}.
\end{proof}

Two edges $(a,c)$ and $(b,d)$ \definition{cross} if $a<b<c<d$.
\begin{lem}\label{L: graph noncrossing}
The edges of an increasing chain graph do not cross.
\end{lem}
\begin{proof}
This follows from Lemma~\ref{L: cycle noncross}.
\end{proof}

\subsection{Two consecutive increasing chains}
We now apply the graph description to the set
$\C(u,w,2) = \{(u,v,w): u \xrightarrow{k} v \xrightarrow{k} w\}$
that is central to this paper. Each triple is represented by the
graphs of its two consecutive increasing chains.
\begin{prop}\label{P: graph pairs}
For $u,w \in S_n$,
the triples $(u,v,w) \in \C(u,w,2)$ are in bijection with graph pairs $(A,B)$ satisfying:
\begin{itemize}
\item Each graph decomposes into ordered paths.
\item We have $\Head(A)=u([k])$ and $\Head(B)=\Tail(A)$.
\item We have $\pi(B)\pi(A)u=w$.
\item We have $|E(A)| + |E(B)| \geq \ell(w) - \ell(u)$.
\end{itemize}

Given such a pair $(A,B)$, we may recover $v$ as $\pi(A)u$.
Then $A$ is the increasing chain graph of $u$ and $v$,
while $B$ is the increasing chain graph of $v$ and $w$.
\end{prop}
\begin{proof}
This follows from Corollary~\ref{C: Graph language}
and the observation that if $v = \pi(A)u$,
then $v([k]) = \Tail(A)$.
\end{proof}

We establish a non-crossing condition
for such $(A,B)$.

\begin{lem}\label{L: mixed crossing reconstruction}
Say $(u,v,w) \in \C(u,w, 2)$ corresponds to $(A,B)$.
There cannot be edges $(a,d)\in E(A)$ and $(b,c)\in E(B)$ with $b<a<c<d$.
\end{lem}
\begin{proof}
Put $i=v^{-1}(a)$ and $j=v^{-1}(c)$, and let
$a'=w(i)$ and $c'=u(j)$.
By Lemma~\ref{L: first graph position values},
either $c = c'$ or $(c, c') \in E(A)$.
In the latter case, Lemma~\ref{L: graph noncrossing} 
guarantees $c' < d$ since $(a,d) \in E(A)$.
In either case, $c' < d$.
Similarly, we can deduce $b < a'$.
In summary,
$$
b < a' \le a < c \le c' < d.
$$
On the other hand, we have
$$
u^{-1}(d)=v^{-1}(a)=w^{-1}(a') \quad \text{ and } \quad
u^{-1}(c')=v^{-1}(c)=w^{-1}(b).
$$
Both positions $i$ and $j$ lie to the right of the cut.
Since $u(j)<u(i)$ and $v(j)>v(i)$,
Proposition~\ref{P:k-Bruhat-endpoints} applied to $u\le_k v$
forces $i<j$.
Similarly, $v(i)<v(j)$ and $w(i)>w(j)$ force $i>j$
by the same proposition applied to $v\le_k w$.
This is a contradiction. 
\end{proof}

\subsection{Increasing chains with double weight}\label{SS:double-weight}

We now consider two sets of variables. 
In addition to the $x$-variables $x_1, \dots, x_m$,
we add $y$-variables $\y = (y_1, \dots, y_n)$.
We introduce a generating function
$F'_{u,w}(x_1,\ldots,x_m;\y)$ 
symmetric in the $x$-variables. Its specialization at $y=0$
is closely related to $F_{u,w}(x_1, \dots, x_m)$
defined earlier.

For $u\xrightarrow{k}w$, put
\[
 \LL_{u,w}=\{u(j):j>k,\ w(j)=u(j)\}.
\]
Since an increasing chain is right-distinct,
\[
 \ell(w)-\ell(u)=n-k-|\LL_{u,w}|.
\]
For $m\ge1$, recall
\[
 \C(u,w,m)=\{(v_0,\ldots,v_m):
 v_0=u,\ v_m=w,\ v_{t-1}\xrightarrow{k}v_t
 \text{ for }1\le t\le m\}.
\]
An element has \definition{double weight}
\[
 \wt_D(v_0, \dots, v_m) := \prod_{t=1}^m\quad\prod_{c \in \LL_{v_{t-1},v_t}} (x_t - y_c).
\]
We define
\begin{equation}\label{EQ: double chain weight}
 F'_{u,w}(x_1,\ldots,x_m;y_1,\ldots,y_n)
 =\sum_{(v_0,\ldots,v_m)\in\C(u,w,m)}
       \wt_D(v_0, \dots, v_m).
\end{equation}

The symmetry follows from Chou--Yu's formulation of Samuel's Pieri identity.
\begin{prop}[{\cite[Theorem 7.1]{Samuel:Pieri},
\cite[Proposition 2.29]{CY:doubleEG}}]\label{P: double weight symmetry}
The polynomial $F'_{u,w}(x_1,\ldots,x_m;\y)$ is symmetric in
$x_1,\ldots,x_m$.
\end{prop}

\begin{lem}\label{L: double weight specialization}
We have
\[
 F'_{u,w}(x_1,\ldots,x_m;\mathbf0)
 =(x_1\cdots x_m)^{n-k}
 F_{u,w}(x_1^{-1},\ldots,x_m^{-1}).
\]
\end{lem}
\begin{proof}
For $(v_0,\ldots,v_m)\in\C(u,w,m)$, put
$r_t=\ell(v_t)-\ell(v_{t-1})$. Since
$|\LL_{v_{t-1},v_t}|=n-k-r_t$, we have
\[
 \prod_{t=1}^m x_t^{n-k-r_t}
 =(x_1\cdots x_m)^{n-k}\prod_{t=1}^m(x_t^{-1})^{r_t}.
\]
Summing over $\C(u,w,m)$ proves the identity.
\end{proof}

\section{Background and preliminaries on positroids}\label{S:positroids}
We recall the matrix and rank descriptions of positroids and record properties of their parallel classes.

\subsection{Positroid varieties, matrices, and bases}
\label{SS: positroid basics}
Fix $u\le_k w$. Let $X_w$ and $X^u$ be the Schubert and opposite Schubert varieties indexed by $w$ and $u$, respectively, in the full flag variety $\Fl(n)$. Their intersection $X_w\cap X^u$ is a \definition{Richardson variety}. Under the projection
\[
 \Fl(n)\longrightarrow\Gr(k,n),\qquad
 (V_1\subseteq\cdots\subseteq V_n)\longmapsto V_k,
\]
its image is the \definition{positroid variety} $\Pi_{u,w}$~\cite{KLS:positroids}.

Choose a full-rank $k\times n$ matrix $X$ whose row span is a general point of $\Pi_{u,w}$. Write $X_a$ for its column in position $a\in[n]$. For $S\subseteq[n]$, put
\[
 r_X(S)=\dim\operatorname{span}\{X_s:s\in S\}.
\]
The resulting function $r_X(\cdot)$ depends only on $u,w$. We write it as $r_{u,w}(S)$, or simply $r(S)$ when $u,w$ are clear.
A set $S\subseteq[n]$ is a \definition{basis} if $r(S)=|S|=k$. The collection of bases, denoted $\mathcal B_{u,w}$, is the \definition{positroid} associated to $\Pi_{u,w}$. An element of $[n]$ is a \definition{loop} if it belongs to no basis, and a \definition{coloop} if it belongs to every basis. An element $c\in[n]$ is a loop exactly when $r(\{c\})=0$, or equivalently when $X_c=0$.

For a $k$-element set $D=\{d_1<\cdots<d_k\}\subseteq[n]$, its \definition{maximal minor} is
\[
 \Delta_D(X)=\det\begin{bmatrix}X_{d_1}&\cdots&X_{d_k}\end{bmatrix}.
\]
Thus
\[
 D\in\mathcal B_{u,w}\quad\Longleftrightarrow\quad\Delta_D(X)\ne0.
\]
A \definition{totally nonnegative realization} of $\mathcal B_{u,w}$ is a real 
$k\times n$ matrix $Y$ whose maximal minors satisfy
\[
 \Delta_D(Y)>0\quad\text{for }D\in\mathcal B_{u,w},\qquad
 \Delta_D(Y)=0\quad\text{for }D\notin\mathcal B_{u,w}.
\]
Such $Y$ exists by~\cite[Proposition 5.11]{KLS:positroids} and~\cite[Theorem 17.1]{Postnikov:total}.

\subsection{Affine permutations and necklaces}
The bounded affine permutation and the Grassmann necklace are two notions that we use to 
compute loops, coloops, and cyclic interval ranks.

For $x\in\mathbb Z$, its \definition{residue} $\overline{x}$ is the unique member of $[n]$ such that $x-\overline{x}$ is divisible by $n$. Given an integer $a$ and $b\in[n]$, there is a unique integer $b' \in [a, a+n)$ with $\overline{b'}=b$.
We call $b'$ the \definition{lift of $b$ in $[a,a+n)$}, or the lift starting at $a$.
Intervals specifying lifts or ranges of integer variables are ordinary integer intervals; intervals of elements of $[n]$ are cyclic as defined below.

A \definition{bounded affine permutation} is a bijection $f:\mathbb Z\to\mathbb Z$ satisfying
\[
 f(a+n)=f(a)+n,\qquad a\le f(a)\le a+n
 \quad\text{for every }a\in\mathbb Z.
\]
Knutson, Lam, and Speyer~\cite[\S3.4]{KLS:positroids} associate the following bounded affine permutation to $u\le_k w$:
\begin{equation}\label{EQ: bounded affine index}
 f_{u,w}(a)=
 \begin{cases}
 u(w^{-1}(a))+n,&\text{if }w^{-1}(a)\le k,\\
 u(w^{-1}(a)),&\text{if }w^{-1}(a)>k,
 \end{cases}
 \qquad a\in[n].
\end{equation}
The identity $f_{u,w}(a+n)=f_{u,w}(a)+n$ determines its remaining values. We write $f$ when $u$ and $w$ are clear. Notice that exactly $k$ numbers among $(f(1),\ldots,f(n))$ exceed $n$.

Extend the column notation by $X_a=X_{\overline a}$ for every $a \in \Z$. The same affine permutation has the description~\cite[\S1.2]{KLS:positroids}
\begin{equation}
\label{EQ: direct f}
 f(a)=\min\{b\ge a:X_a\in\operatorname{span}(X_{a+1},\ldots,X_b)\}.
\end{equation}
The span is zero when $b=a$. Hence $f(a)=a$ exactly when $\overline a$ is a loop. Also, $f(a)=a+n$ exactly when $X_a$ is not in the span of the other $n-1$ columns, which means that $\overline a$ is a coloop.

\begin{lem}\label{L: loops endpoints}
The loops and coloops of $\Pi_{u,w}$ are, respectively,
\[
 \{u(j):j>k,\ u(j)=w(j)\}
 \quad\text{and}\quad
 \{u(i):i\le k,\ u(i)=w(i)\}.
\]
\end{lem}
\begin{proof}
Take $a\in[n]$ and let $q=w^{-1}(a)$. In~\eqref{EQ: bounded affine index}, the equality $f(a)=a$ holds exactly when $q>k$ and $u(q)=w(q)=a$. Similarly, $f(a)=a+n$ holds exactly when $q\le k$ and $u(q)=w(q)=a$.
\end{proof}

Recall that we have defined the set $\LL_{u,w}$ when $u \xrightarrow{k} w$.
We may extend it to $u \le_k w$ by keeping exactly the same definition: 
$ \LL_{u,w} :=\{u(j):j>k,\ u(j)=w(j)\}$.
Now we know this is the set of loops in $\Pi_{u,w}$.

For $a\in[n]$, the \definition{cyclic order starting at $a$} is
\begin{equation}\label{EQ: prec}
 a\prec_a a+1\prec_a\cdots\prec_a n
 \prec_a1\prec_a\cdots\prec_a a-1.
\end{equation}
Equivalently, elements of $[n]$ are compared by their lifts in $[a,a+n)$. The \definition{(Grassmann) necklace} is the sequence of bases $(J_1,\ldots,J_n)$ in which $J_a$ is lexicographically smallest for this order.

For $a,b\in[n]$, the \definition{cyclic interval} $[a,b]$ consists of the numbers encountered from $a$ through $b$ in~\eqref{EQ: prec}, including $a$ and $b$. 
Omitting $a$ or $b$ as appropriate defines $(a,b]$, $[a,b)$, and $(a,b)$. 
In particular, $[a,a]=\{a\}$ and $(a,a]=[a,a)=\varnothing$. We also write 
$r[a,b]=r([a,b])$ for clarity. For integers $a\le b<a+n$, extend this notation by
\[
 r[a,b]=r[\overline a,\overline b].
\]

Put $\rho=wu^{-1}\in S_n$. Equation~\eqref{EQ: bounded affine index} gives, 
for every $x\in\mathbb Z$,
\[
 \overline{f(x)}=\rho^{-1}(\overline x),\qquad
 \overline{f^{-1}(x)}=\rho(\overline x).
\]
We may rewrite Postnikov's necklace correspondence~\cite[\S16, Lemma 16.2]{Postnikov:total} as
\begin{equation}\label{EQ: affine necklace}
 J_a=\{b\in[n]:b\prec_a\rho(b)\}
       \sqcup\{b\in[n]:f(b)=b+n\}.
\end{equation}
The second set consists of the coloops. 
For $i\le k$, we have $\rho(u(i))=w(i)\ge u(i)$. The number $u(i)$ therefore belongs to $J_a$ exactly when $a\notin(u(i),w(i)]$. For $j>k$, $u(j)$ belongs to $J_a$ exactly when $a\in(w(j),u(j)]$. Consequently,
\begin{equation}\label{EQ: positional necklace}
 J_a=\{u(i):i\le k,\ a\notin(u(i),w(i)]\}
 \sqcup\{u(j):j>k,\ a\in(w(j),u(j)]\}.
\end{equation}
Equivalently,
\begin{equation}\label{EQ: positional necklace 2}
 \begin{split}
 J_a={}&\{i \in u([k]): a\notin(i ,\rho(i)]\}\sqcup\{i \in u(\{k+1,\ldots,n\}): a\in(\rho(i),i]\}.
 \end{split}
\end{equation}
In particular, $J_1=u([k])$.

We now recall the following identities from the literature.
\begin{prop}\label{P: positroid ranks}
We have the following statements from the literature.

\begin{enumerate}
\item For $a,b\in[n]$, we have~\cite[Lemma 14]{MO:rank}
\[
 r[a,b]=|J_a\cap[a,b]|.
\]

\item For distinct $a,b\in[n]$, we have~\cite[Theorem 25]{MO:rank}
\begin{equation}\label{EQ: pair rank}
 r(\{a,b\})=
 \min\{r(\{a\})+r(\{b\}),\ r[a,b]+r[b,a]-k\}.
\end{equation}

\item A set $D\subseteq[n]$ of size $k$ is a basis if and only if~\cite[Theorem 6]{Oh:positroids}
\[
 |D\cap[a,b]|\le r[a,b]
 \quad\text{for every cyclic interval }[a,b].
\]
\end{enumerate}
\end{prop}

We record one more way to compute cyclic interval ranks.
\begin{lem}\label{L: affine interval rank}
For integers $a\le b<a+n$,
\begin{equation}\label{EQ: interval rank}
 r[a,b]=
       |\{c\in [a,b]:  f(c) > b\}|.
\end{equation}
If $\overline b$ is not a loop, then equivalently
\begin{equation}\label{EQ: interval rank 2}
 r[a,b]=
      |\{c\in [a,b):  f(c) > b\}| + 1.
\end{equation}
\end{lem}
\begin{proof}
Let $V_{b+1}=\{0\}$, and define 
$V_c=\operatorname{span}(X_c,\ldots,X_b)$ for $c\in[a,b]$.
We study how each $V_c$ differs from $V_{c+1}$:
$V_c = V_{c+1}$ if $X_c \in V_{c+1}$,
or $V_c$ has one more dimension otherwise. 
By~\eqref{EQ: direct f}, $V_c = V_{c+1}$ exactly when $f(c)\le b$. 
Thus, the difference $\dim(V_a) - \dim(V_{b+1})$ 
counts the integers $c\in[a,b]$ for which $f(c)>b$.
Notice that $\dim(V_a) = r[a,b]$ and $\dim(V_{b+1}) = 0$,
which proves~\eqref{EQ: interval rank}. 

If $\overline b$ is not a loop, then $f(b)>b$.
Removing $b$ from the set in~\eqref{EQ: interval rank} 
gives~\eqref{EQ: interval rank 2}.
\end{proof}

\subsection{Parallel classes}
The main notion in this paper is a \definition{parallel class}. As in~\S\ref{sec:introduction}, distinct nonloops $a,b$ are \definition{parallel}, written $a\parallel b$, if $X_a$ and $X_b$ are proportional. Equivalently,
\[
 r(\{a\})=r(\{b\})=r(\{a,b\})=1.
\]
Grouping proportional nonzero columns partitions the elements of $[n]$ that are not loops. Each class is a maximal set whose columns span the same line. Write $\PP_{u,w}$ for the collection of parallel classes.

We record the ``non-crossing'' property of parallel classes.
Four elements $a,b,c,d\in[n]$ are \definition{cyclically ordered} if
\[
 a\prec_a b\prec_a c\prec_a d.
\]

\begin{lem}\label{L: noncrossing parallel classes}
There do not exist cyclically ordered elements $a,b,c,d\in[n]$ in two distinct parallel classes with $a\parallel c$ and $b\parallel d$.
\end{lem}
\begin{proof}
Cyclically rename $a,b,c,d$ so that the smallest is called $a$. Then $a<b<c<d$, $a\parallel c$ and $b\parallel d$.
Choose a totally nonnegative realization $Y$ of $\mathcal B_{u,w}$, and write $\Delta_D=\Delta_D(Y)$. The two classes are distinct, so $\{a,b\}$ is independent and extends to a basis $T\sqcup\{a,b\}$. The set $T$ contains neither $c$ nor $d$, since each would be parallel to a column already in that basis. Write $Tab$ for $T\sqcup\{a,b\}$, and similarly for the other pairs.

Replacing a column by a nonzero parallel column in a basis yields a basis. Consequently $Tab,Tcd,Tad,Tbc$ are bases whereas $\Delta_{Tac}=\Delta_{Tbd}=0$. The determinant identity for these six minors gives
\[
 0=\Delta_{Tac}\Delta_{Tbd}
   =\Delta_{Tab}\Delta_{Tcd}+\Delta_{Tad}\Delta_{Tbc}>0,
\]
a contradiction.
\end{proof}

The next lemma tests whether $a\parallel b$ using only $\rho=wu^{-1}$.

\begin{lem}\label{L: finite parallel test}
For distinct nonloops $a<b$,
\begin{equation}\label{EQ: finite parallel rank}
 r[a,b]+r[b,a]-k
 =1+|\rho((a,b])\setminus[a,b)|.
\end{equation}
Consequently, by~\eqref{EQ: pair rank},
\begin{equation}\label{EQ: parallel transport}
 a\parallel b\quad\Longleftrightarrow\quad\rho((a,b])=[a,b).
\end{equation}
In particular, if $a\parallel b$,
then $a < \rho^{-1}(a) \leq b$.
\end{lem}
\begin{proof}
Applying~\eqref{EQ: interval rank 2} to the integer intervals $[a,b]$ and $[b,a+n]$ gives
\begin{equation}
\label{EQ: r[a,b] r[b,a]}
\begin{split}
r[a,b]&=1 + |\{c\in[a,b):f(c)>b\}|, \\
r[b,a]&=1 + |\{c\in[b,a+n):f(c)>a+n\}|.
\end{split}
\end{equation}

Notice that $r[a,a+n - 1]$ is the rank of the space spanned by all 
columns of $X$, so it is $k$. Applying~\eqref{EQ: interval rank} to $[a,a+n - 1]$ gives
\begin{equation}
\label{EQ: k = c f(c)}
 k=|\{c\in[a,a+n-1]:f(c)>a+n-1\}|.
\end{equation}
The inequalities $x\le f(x)\le x+n$, applied at $x=f^{-1}(a+n)$, give $a\le f^{-1}(a+n)\le a+n$. 
Since $a$ is not a loop, $f^{-1}(a+n) < a+n$, so $f^{-1}(a+n) \in [a,a+n -1]$. Hence,
raising the bound in~\eqref{EQ: k = c f(c)} from $f(c) >  a+n-1$
to $f(c) >  a+n$ decreases the size of the set by $1$. We have
\[
 k=\textcolor{red}{1}+|\{c\in[a,a+n - 1]:  f(c) > \textcolor{red}{a+n}\}|.
\]

Subtracting this equality from the sum in~\eqref{EQ: r[a,b] r[b,a]} yields
\begin{equation}\label{EQ: finite parallel rank sums}
 r[a,b]+r[b,a]-k
  =1+|\{c\in[a,b):f(c)\in(b,a+n]\}|.
\end{equation}

For an integer $c \in [a,b)$, we have $f(c)\le c+n<b+n$. Also, notice that $f(c)>a$: this follows from $f(c)\ge c$ when $c>a$, and from the fact that $a$ is not a loop when $c=a$. Within the range $a<f(c)<b+n$, we have
\[
  f(c)\in(b,a+n]
 \quad\Longleftrightarrow\quad
 \overline{f(c)}\notin(a,b] \quad\Longleftrightarrow\quad \rho^{-1}(c) \notin(a,b].
\]

We rewrite~\eqref{EQ: finite parallel rank sums} as
$$
 r[a,b]+r[b,a]-k
 = 1 + |\{c \in [a,b): \rho^{-1}(c) \notin (a, b]\}|
 = 1 + |\rho((a,b])\setminus[a,b)|
$$
where the last equality uses the bijectivity of $\rho$. 
This proves~\eqref{EQ: finite parallel rank}.

Finally, from~\eqref{EQ: pair rank} we read
$
r(\{a,b\})=\min\{2,\,1+|\rho((a,b])\setminus[a,b)|\}$.
The rank is one exactly when $\rho((a,b])\subseteq[a,b)$. Since these two sets have the same size, $\rho((a,b])\subseteq[a,b)$ is equivalent to 
$\rho((a,b]) = [a,b)$, proving~\eqref{EQ: parallel transport}.
If $a\parallel b$, 
we have $\rho^{-1}(a) \in (a,b]$,
proving the claim at the end. 
\end{proof}

One application is to identify the parallel classes
with connected components of the increasing chain graph 
when $u \xrightarrow{k} w$.

\begin{lem}\label{L: one factor classes}
Suppose $u\xrightarrow{k}w$, with increasing chain graph $G$. Then
$V(G)$ is the set of nonloops, and $\PP_{u,w}$ is the collection of
vertex sets of the connected components of $G$.
\end{lem}

\begin{proof}
By Lemma~\ref{L: first graph position values}, the complement of $V(G)$
is $\LL_{u,w}$, so $V(G)$ is exactly the set of nonloops.
For an edge $(a,b)\in E(G)$, the map $\rho=wu^{-1}=\pi(G)$ sends
$b$ to $a$. Since $a,b$ are consecutive in their component, that
component has no vertex in $(a,b)$. By Lemma~\ref{L: graph noncrossing}, 
if a vertex $v \in V(G)$ is in the interval $(a,b)$,
then vertices in its component all lie in this interval. 
Thus $\rho(v)\in(a,b)\subseteq[a,b)$. Integers outside $V(G)$
are fixed by $\rho=\pi(G)$, so the same inclusion holds for them.
We have $\rho((a,b])\subseteq[a,b)$, so 
$\rho((a,b]) = [a,b)$ and $a\parallel b$ by Lemma~\ref{L: finite parallel test}.

Each component of $G$ is consequently contained in a parallel class.
There are at most $k$ classes, and there are at least $k$ because the
matrix has rank $k$. The components are therefore exactly the 
parallel classes.
\end{proof}

\section{Chain reconstruction algorithm}\label{S: reconstruction}
Fix $u\le_k w$ and put $\rho=wu^{-1}$. Recall that
\[
 I(u,v,w)=\{r:v(u^{-1}(r))=\rho(r)<r\}
 \qquad\text{for }(u,v,w)\in\C(u,w,2).
\]
We define $\Chain(S)$ for a subset $S\subseteq[n]$. The algorithm may
fail on an arbitrary subset. The goal of this section is to show that $\Chain(\cdot)$ is a one-sided inverse of $I(\cdot)$.
\begin{prop}\label{P: reconstruction}
For $(u,v,w)\in\C(u,w,2)$, the algorithm $\Chain$ is well-defined on
$I(u,v,w)$ and gives $(u,v,w)$.
\end{prop}

\subsection*{Describing the algorithm}
We use the identification in Proposition~\ref{P: graph pairs}.
The output of the algorithm is a pair of graphs $(A,B)$. However, we do not assume the output
satisfies the conditions in that proposition.
\begin{defn}\label{D: reconstruct alg}
Take $S\subseteq[n]$. If $S$ contains an $r$ with $\rho(r)\ge r$,
declare $\Chain(S)$ undefined. Otherwise initialize $A_0,B_0$ as empty graphs.
Then the algorithm iteratively 
constructs $A_{r}$ from $A_{r-1}$
and $B_r$ from $B_{r-1}$ for $r = 1, \dots, n$.
During an iteration, the algorithm 
has several ``requirements''. 
If a requirement fails, the algorithm stops and $\Chain(S)$ is undefined.

At iteration $r$, we obtain $(A_r, B_r)$ from $(A_{r-1}, B_{r-1})$ using the first applicable case below.
When adding an edge, also add its endpoints if necessary.
\begin{enumerate}
\item \textbf{Head Case.} If $r\in u([k])$, add the singleton $r$ to $A_{r-1}$ and $B_{r-1}$.
\item \textbf{Null Case.} If $\rho(r)=r$, change neither graph.
\item \textbf{A Case.} If $r\in S$, require that
      \textcolor{red}{$\rho(r)\in\Tail(A_{r-1})$ and $\rho(r)$ is a singleton of $B_{r-1}$}.
      Add $(\rho(r),r)$ to $A_{r-1}$. In $B_{r-1}$, replace the singleton $\rho(r)$ by $r$.
\item \textbf{AB Case.} Otherwise, suppose a path of $A_{r-1}$ has head $h$
      and tail $t$ satisfying
      \[
       \rho(r)<t<r\le\rho(h).
      \]
      Choose the path with minimal $t$. Require that
      \textcolor{red}{$\rho(r)\in\Tail(B_{r-1})$}.
      Add $(t,r)$ to $A_{r-1}$, add $(\rho(r),t)$ to $B_{r-1}$, and add the singleton $r$ to $B_{r-1}$.
\item \textbf{B Case.} Otherwise, require that
      \textcolor{red}{$\rho(r)\in\Tail(B_{r-1})$}, and add $(\rho(r),r)$ to $B_{r-1}$.
\end{enumerate}
At the end of the iteration, we require
\textcolor{red}{$A_{r}$ and $B_{r}$ 
decompose into ordered paths.}
If all iterations succeed, the output is $(A_n,B_n)$. When this pair satisfies
Proposition~\ref{P: graph pairs}, we also identify the output with
$(u,\pi(A_n)u,w)$.
\end{defn}

\begin{exa}\label{Ex: reconstruction scan}
Take $u=[2,1,3,5,4]$, $w=[4,5,1,2,3]$, $k=2$, and $S=\{3\}$.
Iterations $1,2$ are in the Head Case.
The remaining iterations are in the A Case, B Case,
and AB Case, respectively.
We illustrate the graphs as follows.
\[
\begin{array}{c|l|l}
 r&A_r&B_r\\\hline
 0&\varnothing&\varnothing\\
 1&1&1\\
 2&1,\ 2&1,\ 2\\
 3&1-3,\ 2&3,\ 2\\
 4&1-3,\ 2&3-4,\ 2\\
 5&1-3-5,\ 2&2-3-4,\ 5
\end{array}
\]
Notice that $(A_5, B_5)$ satisfy all 
conditions in Proposition~\ref{P: graph pairs},
and correspond to $(u,[2,5,1,3,4],w)$.
\end{exa}

\begin{rem}
\label{R: graph table}
We summarize how each iteration affects
the graphs if it does not fail. We list
how vertices, edges, heads, and tails
of $A_r$ and $B_r$ change from those
of $A_{r-1}$ and $B_{r-1}$.
\[
\begingroup
\renewcommand{\arraystretch}{1.3}
\setlength{\arraycolsep}{4pt}
\begin{array}{c|cccc|cccc}
 &\multicolumn{4}{c|}{A_r}&\multicolumn{4}{c}{B_r}\\
 \text{Case}&V&E&\Head&\Tail&V&E&\Head&\Tail\\\hline
 \text{Head}&+r&-&+r&+r&+r&-&+r&+r\\
 \text{Null}&-&-&-&-&-&-&-&-\\
 \text{A}&+r&+(q,r)&-&-q+r&-q+r&-&-q+r&-q+r\\
 \text{AB}&+r&+(t,r)&-&-t+r&+r&+(q,t)&-t+r&-q+r\\
 \text{B}&-&-&-&-&+r&+(q,r)&-&-q+r
\end{array}
\endgroup
\]
In the table, 
we write $q=\rho(r)$, and $t$ is the chosen tail of $A_{r-1}$
in the AB Case. The notation $+a$ means that $a$ is added,
$-a+b$ means that $a$ is removed and $b$ is added, and a dash means no change.
\end{rem}

\subsection*{Properties of the algorithm}

We start by showing some basic properties
of $(A_r, B_r)$.

\begin{lem}\label{L: reconstruction path components}
If iterations $1,\ldots,r$ succeed, $A_r,B_r$ decompose into ordered paths,
all their vertices belong to $[r]$, and
\[
 \Head(A_r)=u([k])\cap[r],\qquad \Head(B_r)=\Tail(A_r).
\]
\end{lem}
\begin{proof}
Induct on $r$. The base case $r=0$ is immediate.
Now consider the iteration of $r$.
In A Case, the added edge $(\rho(r), r)$
joins $r$ to the ordered path of $\rho(r)$,
so $A_{r}$ decomposes into ordered paths. 
A similar argument shows that $B_r$ 
decomposes into ordered paths in B Case. 
In the AB case, the added edge $(t, r)$
joins $r$ to the ordered path of $t$,
so $A_{r}$ decomposes into ordered paths. 
For $B_r$, recall that 
$t\in\Head(B_{r-1})$ by induction and
$\rho(r)\in\Tail(B_{r-1})$ by the requirement.
Since $\rho(r)<t$, these vertices belong to two distinct paths.
Adding $(\rho(r),t)$ joins these paths into one ordered path, so $B_r$ decomposes into ordered paths. 

All vertices added at this iteration are at most $r$.
The changes in Remark~\ref{R: graph table} give the last
two equations of sets.
\end{proof}

Therefore, the requirement at the end
of each iteration is automatically satisfied.
Next, we show the requirements in AB Case
and B Case are also automatically satisfied. 

\begin{lem}\label{L: available second tail}
Suppose iterations $1,\ldots,r-1$ succeed and $\rho(r)<r$.
Then $\rho(r)\in\Tail(B_{r-1})$.
\end{lem}
\begin{proof}
Iteration $\rho(r)$ is not in the Null Case: otherwise
$\rho(\rho(r))=\rho(r)$, which is impossible.
By Remark~\ref{R: graph table}, this iteration adds $\rho(r)$ to $\Tail(B_{\rho(r)})$.
At iteration $s$ with $\rho(r)<s<r$, the only element that might be removed from $\Tail(B_{s-1})$ is $\rho(s)$.
Thus, $\rho(r)$ remains a tail in 
$B_s$ for $\rho(r)<s<r$.
\end{proof}

Consequently, during an iteration, only the requirements in the A Case might fail. 
We summarize this conclusion in the following corollary.

\begin{cor}
\label{C: requirements checked}
In the algorithm $\Chain$, all requirements during an iteration
are automatically satisfied, except the requirements in A Case.  
\end{cor}
\begin{proof}
Follows from Lemma~\ref{L: reconstruction path components}
and Lemma~\ref{L: available second tail}.
\end{proof}

Notice that the requirements of A Case need not hold for an arbitrary $S$.
\begin{exa}\label{Ex: reconstruction failure}
Take $u=[1,2,3]$, $w=[3,1,2]$ and $k=1$.
Run the algorithm on $S=\{3\}$.
Iteration $1$ adds the singleton $1$ to each graph.
At iteration $2$, $\rho(2)=1$.
We are in the B Case and get: 
\[
 A_2:\quad 1,\qquad B_2:\quad 1-2.
\]
Iteration $3$ is in A Case since $3 \in S$.
However, $\rho(3)=2\notin\Tail(A_2)$.
Thus $\Chain(\{3\})$ is undefined.
\end{exa}

\subsection*{Proof of Proposition~\ref{P: reconstruction}}
In the rest of this section, 
fix $(u,v,w)\in\C(u,w,2)$ and let $A,B$ be the corresponding increasing
chain graphs. Put $S=I(u,v,w)$
and consider running $\Chain(\cdot)$ on $S$.
We give an explicit 
description of 
$V(A_r)$, $E(A_r)$, and $E(B_r)$.

\begin{lem}\label{L: reconstruction restrictions}
When the input is $S$, every iteration succeeds. For $0\le r\le n$,
\begin{equation}\label{EQ: reconstruction restrictions}
\begin{aligned}
 V(A_r)&=V(A)\cap[r],\\
 E(A_r)&=\{(v(j),u(j)):j>k,\ u(j)\le r,\ v(j)<u(j)\},\\
 E(B_r)&=\{(w(j),v(j)):j>k,\ u(j)\le r,\ w(j)<v(j)\}.
\end{aligned}
\end{equation}
\end{lem}
\begin{proof}
Let $V^A_r$, $E^A_r$ and $E^B_r$ be the sets on 
the right hand side of the three equations above respectively.
The goal is to show $V(A_r) = V^A_r$, $E(A_r) = E^A_r$ 
and $E(B_r) = E^B_r$.
We induct on $r$. All sets in the three equations are empty at $r=0$.
Assume the assertion holds through $r-1$
and consider iteration $r$.
The inductive step is immediate if we are in the 
Head Case or Null Case. 

For the remaining cases, we know $r \notin u([k])$
and $\rho(r) < r$.
Let $j=u^{-1}(r)>k$ and $p=v(j)$.
Proposition~\ref{P:k-Bruhat-endpoints} gives $\rho(r)\le p\le r$.
By Lemma~\ref{L: first graph position values},
$r \in V(A)$ if and only if $p < r$. 
By how the sets are defined, 
we have
\begin{equation}\label{EQ: EAEB change}
\begin{aligned}
 \text{ if } p < r, \text{ then} \quad &E^A_r=E^A_{r-1}\sqcup\{(p,r)\},
             \qquad V^A_r=V^A_{r-1}\sqcup\{r\},\\
\text{ if } p = r, \text{ then}\quad &E^A_r=E^A_{r-1},
             \qquad V^A_r=V^A_{r-1},\\
 \text{ if } \rho(r)<p, \text{ then}\quad &E^B_r=E^B_{r-1}\sqcup\{(\rho(r),p)\},\\
 \text{ if } \rho(r)=p, \text{ then} \quad &E^B_r=E^B_{r-1}.
\end{aligned}
\end{equation}
It remains to check $V(A_r)$, $E(A_r)$ and $E(B_r)$
change in the same manner
and the requirements in iteration $r$ are satisfied. 
By Corollary~\ref{C: requirements checked},
we only need to worry about the requirements in A Case. 

\emph{First, suppose $r\in S$,}
so we are in A Case.
Then $\rho(r) = p<r$.
Lemma~\ref{L: first graph position values} says $(p,r)\in E(A)$,
so $p \in V(A)$.
By inductive hypothesis, $p\in V(A_{r-1})$ 
and no edge of the form $(p,s)$ exists in $A_{r-1}$.
Thus, $p\in\Tail(A_{r-1})$.
Lemmas~\ref{L: reconstruction path components} and
\ref{L: available second tail} give
$p\in\Head(B_{r-1})\cap\Tail(B_{r-1})$,
so $p$ is a singleton of $B_{r-1}$.
The requirements of A Case are satisfied. 
By the effect of the algorithm as described in Remark~\ref{R: graph table}, $V(A_r)$, $E(A_r)$ and $E(B_r)$ change from $V(A_{r-1})$, $E(A_{r-1})$ and $E(B_{r-1})$ in the same manner as~\eqref{EQ: EAEB change}.

\emph{Now suppose $r\notin S$,}
so $\rho(r)<p$
which implies $(\rho(r),p)\in E(B)$
by Lemma~\ref{L: first graph position values}.
We claim there is no path in $A_{r-1}$
with head $h$ and tail $t$ 
such that $\rho(r) < t < p$ and $r \leq \rho(h)$.
We prove the claim by considering two cases. 
\begin{itemize}
\item Suppose there is an edge of the form $(t,s)$ in $A$.
Since $t \in \Tail(A_{r-1})$, this edge is not in $A_{r-1}$.
By the inductive hypothesis, $s \geq r$.
If $s = r$, then $(t,r) \in E(A)$ implies $t = v(j) = p$
by Lemma~\ref{L: first graph position values},
contradicting $t < p$.
On the other hand, if $s > r$, we have $\rho(r) < t < p < s$.
Then the edges $(t,s)\in E(A)$ and $(\rho(r),p)\in E(B)$ contradict Lemma~\ref{L: mixed crossing reconstruction}.

\item Suppose no such edge $(t,s)$ exists in $A$,
so $t \in \Tail(A) = \Head(B)$.
The path of $t$ in $B$ cannot have numbers smaller than $t$.
In particular, the path does not have $\rho(r)$.
Since $(\rho(r), p) \in E(B)$,
the path of $t$ does not have $p$ either.
By $\rho(r) < t < p$ and Lemma~\ref{L: graph noncrossing},
the path of $t$ in $B$ lies in the interval $(\rho(r), p)$.

Now recall that $\pi(B)\pi(A)$ is a permutation sending $h$ 
to $\rho(h)$ while $\pi(A)$ sends $h$ to $t$.
Thus, $\pi(B)$ sends $t$ to $\rho(h)$. 
In other words, $\rho(h)$ is the tail of 
the path of $t$ in $B$.
Thus, $\rho(h)$ is in the interval $(\rho(r), p)$,
so $\rho(h)<p\le r$, contradicting $r\le\rho(h)$.
\end{itemize}

This proves our claim. We finish the proof 
by the behavior of the algorithm when $p < r$
and $p = r$.

\emph{Suppose $p<r$.} The edge $(p,r)\in E(A)$ and the induction hypothesis imply $p\in\Tail(A_{r-1})$. Let $h=u(i)$ be the head of this path in $A_{r-1}$,
so $i\le k$. The path of $h$ in $A$ contains $r$ and has tail $v(i)$, so
\[
 \rho(r)<p<r\le v(i)\le w(i)=\rho(h).
\]
Thus this path satisfies the condition of AB Case. Our claim shows that this
is the path satisfying the condition with the smallest tail, 
so the algorithm chooses it. It adds $(p,r)$ to $A_{r-1}$ and $(\rho(r),p)$ to $B_{r-1}$. Thus, $V(A_{r})$, $E(A_{r})$ and $E(B_{r})$ change
in the same manner as~\eqref{EQ: EAEB change}.

Finally, \emph{suppose $p=r$.} Our claim shows that the path
in AB Case cannot exist, so we are in B Case. The iteration adds $(\rho(r),r)=(\rho(r),p)$ to $B_{r-1}$. Again, $V(A_{r})$, $E(A_{r})$ and $E(B_{r})$ change
in the same manner as~\eqref{EQ: EAEB change}.
\end{proof}

\begin{proof}[Proof of Proposition~\ref{P: reconstruction}]
At $r=n$, Lemma~\ref{L: reconstruction restrictions} gives
$V(A_n)=V(A)$, $E(A_n)=E(A)$, and $E(B_n)=E(B)$.
Thus $A_n=A$, and Lemma~\ref{L: reconstruction path components} gives
\[
 \Head(B_n)=\Tail(A_n)=\Tail(A)=\Head(B).
\]
Remark~\ref{R: Two graphs} now gives $B_n=B$.
Therefore $\Chain(I(u,v,w))=(u,v,w)$.
\end{proof}

\begin{cor}
The map $I:\C(u,w,2)\longrightarrow2^{[n]}$ is injective.
\end{cor}
\begin{proof}
This follows from Proposition~\ref{P: reconstruction}.
\end{proof}

\section{Length bounds and increasing chains}\label{S: inequalities}

We prove two lower bounds for $\ell(w)-\ell(u)$. Equality in the
first characterizes the existence of one increasing chain,
$u\xrightarrow{k}w$. Equality in the second guarantees a factorization
$u\xrightarrow{k}v\xrightarrow{k}w$ into two increasing chains.

For $u,w\in S_n$, define
\[
 m(u,w)=|\{j\in\{k+1,\ldots,n\}:u(j)>w(j)\}|.
\]
When $u\le_k w$, Proposition~\ref{P:k-Bruhat-endpoints} gives
\[
 m(u,w)=|\{j\in[n]:u(j)>w(j)\}|
       =|\{r\in[n]:r>\rho(r)\}|,
 \qquad \text{ where } \rho=wu^{-1}.
\]

Our first bound assumes only certain inequalities between entries.
In particular, it does not assume $u\le_k w$.

\begin{prop}\label{P: k-Bruhat easy check}
Let $u,w\in S_n$ satisfy
\[
 u(i)\le w(i)\quad(i\le k),\qquad
 u(j)\ge w(j)\quad(j>k).
\]
Then
$\ell(w)-\ell(u)\ge m(u,w)$
with equality if and only if $u\xrightarrow{k}w$.
\end{prop}

The condition $u\xrightarrow{k}w$ is equivalent to
$\C(u,w,1)\ne\varnothing$.
For $u\le_k w$, we use parallel classes 
to give a stronger bound.
We will show equality in this bound is equivalent 
to $\C(u,w,2) \neq \varnothing$.

\begin{prop}\label{P: numerical class bound}\label{L: numerical class bound}
For $u\le_k w$,
\begin{equation}\label{EQ: numerical class bound}
 \ell(w)-\ell(u)\ge m(u,w)+|\PP_{u,w}|-k.
\end{equation}
If equality holds, then $\C(u,w,2)\ne\varnothing$.
\end{prop}

The converse of the last assertion is proved in
Corollary~\ref{C: two factor equality} of the next section.
The new bound is at least as strong as the first
since $|\PP_{u,w}|\ge k$.

We first prove the weaker bound. We then recall Hall's theorem and
derive the form used in the proof of
Proposition~\ref{P: numerical class bound}.

\subsection{The weaker bound}\label{SS: one chain bound}
We begin with a basic lemma on increasing chains. 
\begin{lem}\label{L: parabolic chain transport}
Let $s=t_{i,i+1}$ with $1\le i<n$ and $i\ne k$. Suppose
\[
 \ell(ws)-\ell(us)=\ell(w)-\ell(u).
\]
Then $u\xrightarrow{k}w$ if and only if $us\xrightarrow{k}ws$.
\end{lem}
\begin{proof}
Suppose $v_0\lessdot_k\cdots\lessdot_k v_m$ is the increasing chain
from $u$ to $w$, with $v_r=t_{a_r,b_r}v_{r-1}$ and
$a_1<\cdots<a_m$. 
We show that $(v_0s, v_1s, \dots, v_ms)$
is an increasing chain from $us$ to $ws$.
Clearly, $v_rs=t_{a_r,b_r}v_{r-1}s$.
Since $s = t_{i, i+1}$ with $i \neq k$,
$t_{a_r,b_r}$ is a $k$-swap on $v_{r-1}s$.
It remains to show that it is a $k$-cover.

Lemma~\ref{L:cover-test} gives
$\ell(v_rs)-\ell(v_{r-1}s)\ge1$ for every $r$.
The sum of these $m$ increments is
\[
\ell(v_ms) - \ell(v_0s) = \ell(ws) - \ell(us)
= \ell(w) - \ell(u) = m
\]
by assumption. This forces
$\ell(v_rs)-\ell(v_{r-1}s) = 1$,
so $t_{a_r,b_r}$ is a $k$-cover on $v_{r-1}s$.
Consequently, $(v_0s, v_1s, \dots, v_ms)$
is an increasing chain from $us$ to $ws$.
The converse follows by applying the same argument to $us$ and $ws$,
using $s^2=1$.
\end{proof}

\begin{proof}[Proof of Proposition~\ref{P: k-Bruhat easy check}]
Let $u'\in S_n$ be obtained by arranging the entries of $u$
in decreasing order within positions $1,\ldots,k$ and within positions
$k+1,\ldots,n$. 
We can find a saturated chain
\[
 u=u_0\lessdot u_1\lessdot\cdots\lessdot u_p=u'
\]
such that each $u_i$ is obtained from $u_{i-1}$ by exchanging the entries in two adjacent
positions on the same side of the cut. 
Say $u_i = u_{i-1}t_{j_i, j_i + 1}$, with $j_i\ne k$.
Now consider the sequence $w_0, w_1, \dots, w_p$
where $w_0 = w$ and $w_i = w_{i-1}t_{j_i, j_i + 1}$.
Thus, for $i \in [p]$,
$\ell(w_i) = \ell(w_{i-1}) + 1$ or $\ell(w_{i-1}) - 1$,
but $\ell(u_i) = \ell(u_{i-1}) + 1$.
Let $w' = w_p$, so
\[
 \ell(w)-\ell(u)\ge \ell(w')-\ell(u').
\]

The coordinate inequalities and $m(u,w)$ are preserved by the common swaps, so
\[
 u'(i)\le w'(i)\text{ for }i\le k \text{ and }
 u'(j)\ge w'(j)\text{ for }j>k.
\]
Since $u'$ decreases within each position block, the within-block
condition in Proposition~\ref{P:k-Bruhat-endpoints} is vacuous. Consequently $u'\le_k w'$. In any saturated chain
from $u'$ to $w'$, each position $j$ with $u'(j) > w'(j)$ 
must be involved in a swap, so 
\[
 \ell(w)-\ell(u)\ge \ell(w')-\ell(u') \ge |\{j: u'(j) > w'(j)\}|
 = |\{j: u(j) > w(j)\}|.
\]

Finally, suppose $ \ell(w)-\ell(u) = |\{j: u(j) > w(j)\}|$,
so the inequalities above are all equalities.
Denote this number by $m$.
Take a saturated chain in $k$-Bruhat order from $u'$ to $w'$,
so it has $m$ swaps.
Each of the $m$ positions $j$ with $u'(j)>w'(j)$ must be involved
in at least one swap.
No position greater than $k$ can occur twice. Thus the chain is
right-distinct, and hence $u'\xrightarrow{k}w'$.

In addition, $\ell(w)-\ell(u)=\ell(w')-\ell(u')$ implies
$\ell(w_i)=\ell(w_{i-1})+1$ for each $i\in[p]$.
Applying Lemma~\ref{L: parabolic chain transport} successively
gives $u_i \xrightarrow{k} w_i$ for $i = p, \dots, 1,0$.
In particular $u \xrightarrow{k} w$.

Conversely, an increasing chain is right-distinct. Each position greater
than $k$ that changes from $u$ to $w$ occurs exactly once, and no other such
position occurs. Its length is therefore $m(u,w)$.
This proves the converse.
\end{proof}

\subsection{Hall's marriage theorem}\label{SS: coordinate interpolation}
To prove Proposition~\ref{P: numerical class bound}, we will construct a
permutation by choosing one integer from each of several sets, with no
integer chosen twice. Hall's theorem characterizes when this is possible.

\begin{thm}[{\cite{Hall:representatives}}]
\label{T: Hall representatives}
Let $Q_1,\ldots,Q_s\subseteq[s]$. There are distinct integers
$q_1,\ldots,q_s$ with $q_i\in Q_i$ if and only if
\begin{equation}\label{EQ: Hall 1}
 |J|\le\left|\bigcup_{j\in J}Q_j\right|
 \qquad\text{for every }J\subseteq[s].
\end{equation}
\end{thm}

We derive the following variations.
\begin{lem}\label{L: Hall contained sets}
Condition~\eqref{EQ: Hall 1} can be replaced by
\begin{equation}\label{EQ: Hall 2}
 |\{i\in[s]:Q_i\subseteq S\}|\le |S|
 \qquad\text{for every }S\subseteq[s].
\end{equation}
If each $Q_i$ is an interval, with $\varnothing$ also considered an interval,
then condition~\eqref{EQ: Hall 2} is equivalent to
\begin{equation}\label{EQ: Hall 3}
 |\{i\in[s]:Q_i\subseteq S\}|\le |S|
 \qquad\text{for every interval }S\subseteq[s].
\end{equation}
\end{lem}
\begin{proof}
Assume~\eqref{EQ: Hall 1}. Given $S\subseteq[s]$, put
$J=\{i:Q_i\subseteq S\}$. Then
\[
 |\{i:Q_i\subseteq S\}|=|J|
 \le\left|\bigcup_{i\in J}Q_i\right|\le |S|,
\]
which proves~\eqref{EQ: Hall 2}.
Conversely, given $J\subseteq[s]$, put $S=\bigcup_{j\in J}Q_j$.
Every $j\in J$ satisfies $Q_j\subseteq S$, so~\eqref{EQ: Hall 2} gives
\[
 |J|\le |\{i:Q_i\subseteq S\}|\le |S|
 =\left|\bigcup_{j\in J}Q_j\right|.
\]

Now suppose each $Q_i$ is an interval.
Condition~\eqref{EQ: Hall 2} is clearly stronger 
than~\eqref{EQ: Hall 3}.
Conversely, suppose~\eqref{EQ: Hall 3} holds. Taking $S=\varnothing$
shows that every $Q_i$ is nonempty. Given a nonempty $S\subseteq[s]$,
write
\[
 S=S_1\sqcup\cdots\sqcup S_t
\]
as the disjoint union of maximal intervals contained in $S$. 
Each $Q_i$ contained in $S$ lies in exactly one $S_j$. Thus
\[
 |\{i:Q_i\subseteq S\}|
 =\sum_{j=1}^t|\{i:Q_i\subseteq S_j\}|
 \le\sum_{j=1}^t|S_j|=|S|.
\]
This proves~\eqref{EQ: Hall 2}.
\end{proof}

Finally, we generalize the statement from $[s]$ to $M \subseteq[n]$.
For such $M$, an \definition{interval subset of $M$} is the empty
set or a subset consisting of consecutive members of $M$. 
Every nonempty such subset has the form $M\cap[a,b]$ for integers $a\le b$.

\begin{cor}\label{C: Hall}
Let $M\subseteq[n]$, and let $Q_1,\ldots,Q_{|M|}$ be interval subsets of $M$.
Suppose
\begin{equation}\label{EQ: interval assignment capacity}
 |\{i:Q_i\subseteq S\}|\le |S|
 \qquad\text{for every interval subset }S\text{ of }M.
\end{equation}
 Then there are distinct integers
$q_1,\ldots,q_{|M|}$ with $q_i\in Q_i$ for every $i\in[|M|]$.
\end{cor}
\begin{proof}
List $M=\{a_1 < \cdots < a_{|M|}\}$ and replace each $a_j$ by $j$.
This gives a bijection from $M$ to $[|M|]$.
An interval subset of $M$ is identified with an interval
of $[|M|]$ with the same size.
The result follows from using condition~\eqref{EQ: Hall 3}
in Theorem~\ref{T: Hall representatives}.

\end{proof}

\subsection{The stronger bound}\label{SS: numerical preliminaries}
Fix $u\le_k w$, put $\rho=wu^{-1}$, and let
\[
 M=\{\min C:C\in\PP_{u,w}\},
\]
so $|M|=|\PP_{u,w}|$. We first record two properties of $M$.
\begin{lem}\label{L: minima interval escape}
For $a\in M\cup\{0\}$ and $b\in M$ with $a<b$,
\[
 \rho(M\cap(a,b])\not\subseteq[a,b).
\]
\end{lem}
\begin{proof}
Suppose the lemma fails. 
Choose a counterexample $a,b$ with $a$ as small as possible,
so $\rho(M\cap(a,b])\subseteq[a,b)$.
We claim that $\rho((a,b])\subseteq[a,b)$.
For $r\in(a,b]$, consider the following three cases.
\begin{itemize}
\item If $r\in M$, then $\rho(r)\in[a,b)$ by assumption.
\item If $r$ is a loop, then $\rho(r)=r$. Since $b\in M$ is not a loop,
we have $r<b$, and hence $\rho(r)\in[a,b)$.
\item Otherwise, let $r'<r$ be the minimum of the class containing $r$.
Lemma~\ref{L: finite parallel test} gives
$\rho((r',r])=[r',r)$.
If $r'<a$, then
\[
 \begin{aligned}
 \rho(M\cap(r',b])
 &\subseteq \rho((r',r])\cup\rho(M\cap(a,b])\\
 &\subseteq [r',r)\cup[a,b)=[r',b).
 \end{aligned}
\]
Thus $r',b$ contradict the minimality of $a$.
Consequently $r'\ge a$, and
$\rho(r)\in[r',r)\subseteq[a,b)$.
\end{itemize}
This proves the claim. The two sets $\rho((a,b])$
and $[a,b)$
have the same
size, so they 
are equal.
If $a=0$, $0$ is in the image of $\rho$,
which is impossible.
If $a>0$, Lemma~\ref{L: finite parallel test} gives $a\parallel b$,
contradicting the fact that $b$ is the minimum of its class.
\end{proof}

\begin{lem}\label{L: heads in M}
We have $u([k])\subseteq M$.
\end{lem}
\begin{proof}
Take $b\in u([k])$. By Lemma~\ref{L: loops endpoints}, $b$ is a nonloop.
If $a<b$ belonged to the same class, Lemma~\ref{L: finite parallel test}
would give $\rho((a,b])=[a,b)$ and hence $\rho(b)<b$.
However, $u^{-1}(b)\le k$ implies
$\rho(b)=w(u^{-1}(b))\ge u(u^{-1}(b))=b$.
Thus $b$ is the minimum of its class, so $b\in M$.
\end{proof}

For $r\in u([k])$, we have
$\rho(r)\ge r$.
For $r\in M\setminus u([k])$, we have $u^{-1}(r)>k$ and $r$ is a
nonloop. Hence $\rho(r)<r$ by Lemma~\ref{L: loops endpoints}.
For each $r\in M$, define
\[
 Q_r=\begin{cases}
 M\cap[r,\rho(r)],&r\in u([k]),\\
 M\cap(\rho(r),r),&r\in M\setminus u([k]).
 \end{cases}
\]
Each $Q_r$ is an interval subset of $M$
in the sense of 
Corollary~\ref{C: Hall}. The following lemma verifies
the condition in Corollary~\ref{C: Hall}.
\begin{lem}
\label{L: Hall on Q_r}
Let $a\in M\cup\{0\}$ and $b\in M\cup\{n+1\}$ with $a<b$,
and put $S=M\cap(a,b)$. For every $r\in M$,
\begin{equation}\label{EQ: minima Hall block}
 Q_r\subseteq S
 \quad\Longleftrightarrow\quad
 r\in (a,b]\ \text{and}\ \rho(r)\in[a,b).
\end{equation}

Consequently, $|\{r\in M:Q_r\subseteq S\}|\le |S|$.
\end{lem}
\begin{proof}
We first prove~\eqref{EQ: minima Hall block} in two cases.
\begin{itemize}
\item Suppose $r\in u([k])$. Then $\rho(r)\ge r$ and
$Q_r=M\cap[r,\rho(r)]$, whose smallest member is $r$.
If $Q_r\subseteq S$, then $a<r<b$. When $b\in M$, the fact that
$b\notin Q_r$ gives $\rho(r)<b$. 
When $b = n+1$, $\rho(r) < b$ is immediate. 
Thus $a<r\le\rho(r)<b$,
which gives the right-hand side of~\eqref{EQ: minima Hall block}.

Conversely, $r\in(a,b]$ and $\rho(r)\in[a,b)$, together with
$r\le\rho(r)$, imply $[r,\rho(r)]\subseteq(a,b)$, so $Q_r\subseteq S$.

\item Suppose $r\in M\setminus u([k])$. Then $\rho(r)<r$ and
$Q_r=M\cap(\rho(r),r)$. Let $r^-$ be the largest member of
$M\cup\{0\}$ smaller than $r$. Since $M\cap(r^-,r]=\{r\}$,
Lemma~\ref{L: minima interval escape} gives
$\rho(r)\notin[r^-,r)$. Therefore $\rho(r)<r^-<r$.
In particular, $r^->0$, so $r^-\in M$ and it is the greatest member of $Q_r$.

Suppose $Q_r\subseteq S$. 
Then $r^- \in Q_r \subseteq S$ implies $a<r^-<b$.
We claim $r \le b$: Otherwise, $r^- < b < r$ and $b \in M$
contradict our choice of $r^-$.
We have $a < r^- < r \leq b$.
It remains to check $\rho(r)\in[a,b)$.
We already have $\rho(r) < r \leq b$.
If $a\in M$, then $a \notin S$ implies 
$a \notin Q_r$, so $a \leq \rho(r)$.
If $a = 0$, $a \leq \rho(r)$ is immediate. 

Conversely, suppose $a<r\le b$ and $\rho(r)\in[a,b)$.
We have $(\rho(r), r) \subseteq (a,b)$,
so $Q_r \subseteq S$.
\end{itemize}

We now count the integers $r\in M$ satisfying $Q_r\subseteq S$.
By~\eqref{EQ: minima Hall block},
such $r$ must come from $M \cap (a,b]$,
which is $S$ when $b = n+1$ or $S \sqcup \{b\}$
when $b \in M$.
Suppose, for a contradiction, that more than $|S|$ integers $r\in M$
satisfy $Q_r\subseteq S$. Then $b\in M$, and every member of
$M\cap(a,b]$ satisfies this containment.
By~\eqref{EQ: minima Hall block}, for 
every $r \in M \cap (a,b]$, 
$\rho(r) \in [a,b)$. 
This contradicts Lemma~\ref{L: minima interval escape}.
\end{proof}

Every interval subset of $M$, including the empty set, has the form
$M\cap(a,b)$ used in Lemma~\ref{L: Hall on Q_r}.
We can therefore apply Corollary~\ref{C: Hall} to construct the following permutation.

\begin{lem}\label{L: minima midpoint}\label{L: coordinate interpolation}
There is a permutation $x\in S_n$ such that
\begin{equation}\label{EQ: u x w}
\begin{aligned}
 u(i)&\le x(i)\le w(i)&& (i\le k),\\
 u(j)&\ge x(j)\ge w(j)&& (j>k),
 \end{aligned}
\end{equation}
and
\[
 m(u,x)=|M|-k,\qquad m(x,w)=m(u,w).
\]
\end{lem}

\begin{proof}
By Lemma~\ref{L: Hall on Q_r} and Corollary~\ref{C: Hall}, choose
distinct integers $q_r\in Q_r$ for $r\in M$.
They use each member of $M$ exactly once, so
\[
 x(u^{-1}(r))=\begin{cases}
 q_r,&r\in M,\\
 r,&r\notin M
 \end{cases}
\]
defines a permutation. The definitions of $Q_r$ give~\eqref{EQ: u x w}. For $j>k$,
$u(j)>x(j)$ exactly when $j = u^{-1}(r)$ 
for $r\in M\setminus u([k])$.
Thus, $m(u,x)=|M|-k$. 
At each such position,
$x(j) > w(j)$ by construction. 
Hence $x(j)>w(j)$ if and only if
$u(j)>w(j)$ for all $j>k$, proving $m(x,w)=m(u,w)$.
\end{proof}

\begin{exa}\label{Ex: concrete Hall midpoint}
Continue Example~\ref{Ex: three parallel pairs},
we take $u=[3,1,6,5,4,2]$, $w=[6,4,5,2,3,1]$, and $k=2$.
We compute:
\[
 M=\{1,3,5\},\qquad u([2])=\{1,3\},\qquad
 \begin{array}{c|ccc}
 r&1&3&5\\\hline
 \rho(r)&4&6&2\\
 Q_r&\{1,3\}&\{3,5\}&\{3\}
 \end{array}.
\]
We may pick $q_1=1$, $q_3=5$ and $q_5=3$.
Notice that they are three distinct numbers
and $q_r \in Q_r$.
The construction in the proof of Lemma~\ref{L: minima midpoint} gives
\[
 x=[5,1,6,3,4,2].
\]
We have $m(u,x)=1=|M|-k$ and $m(x,w)=4=m(u,w)$.
One may verify that $x$ satisfies the conditions
in Lemma~\ref{L: minima midpoint}.
\end{exa}

\begin{proof}[Proof of Proposition~\ref{P: numerical class bound}]
Choose $x$ as in Lemma~\ref{L: minima midpoint}.
By~\eqref{EQ: u x w},
we can apply 
Proposition~\ref{P: k-Bruhat easy check} to $u,x$ and to $x,w$,
obtaining
\[
 \ell(x)-\ell(u)\ge m(u,x) =  |M|-k,\qquad
 \ell(w)-\ell(x)\ge m(x,w) = m(u,w).
\]
Adding them proves~\eqref{EQ: numerical class bound}
since $|M|=|\PP_{u,w}|$.

If equality holds in~\eqref{EQ: numerical class bound}, both
inequalities above are equalities. Proposition~\ref{P: k-Bruhat easy check} then gives
$u\xrightarrow{k}x\xrightarrow{k}w$, so $\C(u,w,2)\ne\varnothing$.
\end{proof}

\section{Changes in parallel classes}\label{S: parallel classes}
We prove that a factorization into two increasing chains forces equality
in Proposition~\ref{P: numerical class bound}. The main step describes
how the parallel classes change during the second increasing chain.

Fix $(u,v,w)\in\C(u,w,2)$, with increasing chain graphs $A,B$.
Write the chain from $v$ to $w$ as
\begin{equation}\label{EQ: chain from v to w}
 v=z_0\lessdot_k z_1\lessdot_k\cdots\lessdot_k z_p=w,
 \qquad z_s=t_{i_s,j_s}z_{s-1},
 \qquad i_1<\cdots<i_p,
\end{equation}
where $i_s<j_s$. For $1\le s\le p$, put
\[
 i'_s=u(z_{s-1}^{-1}(i_s)),\qquad
 j'_s=u(z_{s-1}^{-1}(j_s)).
\]

Lemma~\ref{L: one factor classes} describes $\PP_{u,v}$.
The following proposition describes how each $\PP_{u,z_s}$
differs from $\PP_{u,z_{s-1}}$.

\begin{prop}\label{P: second factor class update}
For $1\le s\le p$, the classes of $\PP_{u,z_s}$ are obtained from
those of $\PP_{u,z_{s-1}}$ as follows.
\begin{enumerate}
\item If $j_s\notin V(A)$, 
then $j_s$ is adjoined to one existing class
as its largest member.
\item If $j_s\in V(A)$, then $j_s,j'_s$ belong to one class
$C\in\PP_{u,z_{s-1}}$, which is replaced by
$C\cap(i_s,j_s]$ and $C\setminus(i_s,j_s]$.
\end{enumerate}
Every other class is unchanged. 
Moreover, for $0\le s\le p$,
\begin{equation}\label{EQ: two factor class equality}
 \ell(z_s)-\ell(u)=m(u,z_s)+|\PP_{u,z_s}|-k.
\end{equation}
\end{prop}

Together with Proposition~\ref{P: numerical class bound}, this gives the
following characterization of $\C(u,w,2)\ne\varnothing$.

\begin{cor}\label{C: two factor equality}
For $u\le_k w$,
\[
 \C(u,w,2)\ne\varnothing
 \quad\Longleftrightarrow\quad
 \ell(w)-\ell(u)=m(u,w)+|\PP_{u,w}|-k.
\]
\end{cor}
\begin{proof}
The forward implication is~\eqref{EQ: two factor class equality} at
$s=p$. The converse is
Proposition~\ref{P: numerical class bound}.
\end{proof}

We first compare the parallel classes for an arbitrary cover
$u\le_k z\lessdot_k z'$. We then apply that comparison to
\eqref{EQ: chain from v to w}.

\subsection{Parallel classes under one \texorpdfstring{$k$}{k}-Bruhat cover}\label{SS: general class update}
Throughout this subsection, assume
\[
 u\le_k z\lessdot_k z'=t_{i,j}z,\qquad i<j.
\]
Let
\[
\sigma=zu^{-1},\qquad \sigma'=z'u^{-1},\qquad
 i_0=\sigma^{-1}(i),\qquad j_0=\sigma^{-1}(j).
\]
The cover uses positions $z^{-1}(i)\le k<z^{-1}(j)$.
The coordinate comparisons for $u\le_k z$ give
$i_0\le i<j\le j_0$.
We will repeatedly use the following two cyclic intervals: 
\[
 (i,j]=\{i+1,\ldots,j\},\qquad
 [j_0,i_0)=\{j_0,\ldots,n\}\cup\{1,\ldots,i_0-1\}.
\]
For $a\in(i,j]$, we have $j_0\prec_a i_0$.

Let $X$ be a $k\times n$ matrix whose row span is a general point of $\Pi_{u,z}$. 
In this subsection, put $f=f_{u,z}$ and write $r,r'$ for the rank functions
of $\Pi_{u,z},\Pi_{u,z'}$, respectively. 
We compare the loops and parallel classes of these two positroid
varieties.

\begin{lem}
\label{L: loop change}
We have
\[
 \LL_{u,z'}=
 \begin{cases}
 \LL_{u,z},&j<j_0,\\
 \LL_{u,z}\setminus\{j\},&j=j_0.
 \end{cases}
\]
In particular, every nonloop of $\Pi_{u,z}$ remains
a nonloop in $\Pi_{u,z'}$.
\end{lem}
\begin{proof}
For $r\in u(\{k+1,\ldots,n\})$, the values $\sigma(r)$ and $\sigma'(r)$
differ only at $r=j_0$. Also,
\[
 \sigma(j_0)=j,\qquad \sigma'(j_0)=i<j\le j_0.
\]
If $j=j_0$, then $j_0=j$ belongs to $\LL_{u,z}$ but not to $\LL_{u,z'}$.
If $j<j_0$, then $j_0$ belongs to neither set. All other memberships agree.
\end{proof}

We next compare $\PP=\PP_{u,z}$ and $\PP'=\PP_{u,z'}$.

\begin{prop}\label{P: general cover parallel classes}
The collections $\PP$ and $\PP'$ are related as follows.
\begin{enumerate}
\item Suppose $j<j_0$. There is at most one $C\in\PP$ meeting both $(i,j]$ and $[j_0,i_0)$. If there is no such class, then $\PP'=\PP$. Otherwise $C\subseteq (i,j]\cup [j_0,i_0)$ and
\[
 \PP'=(\PP\setminus\{C\})\cup\{C\cap (i,j],C\cap [j_0,i_0)\}.
\]
Thus this class splits into two classes.

\item Suppose $j=j_0$. No two existing classes merge, and no existing class splits. Define
\[
 Q=\{a\in[n]:r[a,a]=1,\ i_0\le a\le i,\ r[a,j]+r[j,a]=k+1\}.
\]
If $Q=\varnothing$, then $\PP'=\PP\cup\{\{j\}\}$. 
Otherwise $Q$ is one class in $\PP$ and
\[
 \PP'=(\PP\setminus\{Q\})\cup\{Q\cup\{j\}\}.
\]
In particular, $j$ is the largest element of its class in $\PP'$.
\end{enumerate}
\end{prop}

\begin{exa}\label{Ex: general cover split other class}
Take $u=[3,2,1,6,5,4]$, $z=[5,6,1,4,3,2]$, and $k=3$. Let $z'=t_{1,3}z=[5,6,3,4,1,2]$, so $i_0=1$ and $j_0=5$. The class $\{2,6\}\in\PP$ splits into $\{2\}$ and $\{6\}$ in $\PP'$, while the other classes remain unchanged.
\end{exa}

\begin{exa}\label{Ex: general cover joins class}
Take $k=2$, $u=[1,2,3,4]$, and $z=[1,3,2,4]$. Let$z'=t_{3,4}z=[1,4,2,3]$.
Here $i=3$, $i_0=2$, and $j_0=j=4$. The classes 
in $\PP$ are $\{1\}$ and $\{2,3\}$,
while $4$ is a loop. In $\PP'$ 
the classes are $\{1\}$ and $\{2,3,4\}$.
In other words, $4$ used to be a loop and 
now joins the class $Q=\{2,3\}$ as its largest element.
\end{exa}

\begin{exa}\label{Ex: general cover new singleton}
Take $u=[3,2,1,5,4]$, $z=[4,3,2,5,1]$, and $k=3$. Let $z'=t_{3,5}z=[4,5,2,3,1]$, so $j_0=j=5$. The classes of $\PP$ are $\{1\},\{2\},\{3\}$ and $\{4\}$. Then $\PP'$ is obtained by adding the new singleton class $\{5\}$.
\end{exa}

We now prove Proposition~\ref{P: general cover parallel classes}.
Let $(J_1,\ldots,J_n)$ and $(J'_1,\ldots,J'_n)$ be the necklaces of $\Pi_{u,z}$ and $\Pi_{u,z'}$. We first compare the two necklaces.
\begin{lem}\label{L: one cover necklace update}
We have
\[
 J'_a=
 \begin{cases}
 (J_a\setminus\{i_0\})\cup\{j_0\},&\text{if }i<a\le j,\\
 J_a,&\text{otherwise}.
 \end{cases}
\]
\end{lem}
\begin{proof}
Use the two sets in~\eqref{EQ: positional necklace 2}. For
$x\in u([k])$, the interval $(x,\sigma(x)]$ changes only when $x=i_0$:
it changes from $(i_0,i]$ to $(i_0,j]$. The first set therefore loses
$i_0$ exactly when $a\in(i,j]$.
For $x\in u(\{k+1,\ldots,n\})$, the interval $(\sigma(x),x]$ changes
only when $x=j_0$: it changes from $(j,j_0]$ to $(i,j_0]$.
The second set therefore gains $j_0$ exactly when $a\in(i,j]$.
\end{proof}

\begin{lem}\label{L: cover rank rectangle}
\label{C: cover separation criterion}\label{L: cover separation criterion}
We have the following comparisons.
\begin{enumerate}
\item For $a,b\in[n]$,
\begin{equation}\label{EQ: class update rank}
 r'[a,b]-r[a,b]=\mathbf1_{\{a\in (i,j]\}}\mathbf1_{\{b\in [j_0,i_0)\}}.
\end{equation}
Here $\mathbf1_{\mathcal A}$ is $1$ if $\mathcal A$ holds and $0$
otherwise.
\item If $a,b$ are in two different classes of $\PP$,
they remain in different classes in $\PP'$. 
If $a,b$ are in one class of $\PP$,
then they belong to different classes in $\PP'$ exactly when one belongs
to $(i,j]$ and the other to $[j_0,i_0)$.
\item If $j<j_0$ and a class $C\in\PP$ meets both $(i,j]$ and $[j_0,i_0)$,
then $C\subseteq (i,j]\cup [j_0,i_0)$. Moreover, $C\cap (i,j]$ and $C\cap [j_0,i_0)$ are classes
in $\PP'$. Any class not meeting both sets remains a class in $\PP'$.
\end{enumerate}
\end{lem}
\begin{proof}
For (1), Proposition~\ref{P: positroid ranks}(1) gives
$r[a,b]=|J_a\cap[a,b]|$ and $r'[a,b]=|J'_a\cap[a,b]|$.
By Lemma~\ref{L: one cover necklace update}, these cardinalities agree
unless $a\in (i,j]$. In that case $J'_a=(J_a\setminus\{i_0\})\cup\{j_0\}$.
In the cyclic order starting at $a$, the number $j_0$ occurs before $i_0$.
Consequently $[a,b]$ contains $j_0$ but not $i_0$ exactly when $b\in [j_0,i_0)$.
This proves~(1).

Now suppose $a,b$ are distinct nonloops of $\Pi_{u,z}$,
so they remain nonloops of $\Pi_{u,z'}$ by Lemma~\ref{L: loop change}.
By~\eqref{EQ: pair rank}, the sum $r[a,b]+r[b,a]-k$ equals $1$
exactly when $a,b$ are parallel, and is greater than $1$ otherwise.
Part~(1) gives
\[
\begin{aligned}
 r'[a,b]+r'[b,a]-k
 ={}&r[a,b]+r[b,a]-k\\
 &+\mathbf1_{\{a\in(i,j]\}}\mathbf1_{\{b\in[j_0,i_0)\}}\\
 &+\mathbf1_{\{b\in(i,j]\}}\mathbf1_{\{a\in[j_0,i_0)\}}.
\end{aligned}
\]

The two added terms are nonnegative. Therefore $a,b$ remain
parallel exactly when they were parallel and both added terms vanish.
They are separated exactly when one lies in $(i,j]$ and the other in
$[j_0,i_0)$. This proves~(2).

Finally, suppose $j<j_0$, so $(i,j]\cap[j_0,i_0)=\varnothing$.
Take a class $C\in\PP$ containing $a\in(i,j]$ and $b\in[j_0,i_0)$.
If $c\in C\setminus((i,j]\cup[j_0,i_0))$, then part~(2) would make
$c$ parallel to both $a$ and $b$ in $\PP'$, while separating $a$
from $b$. This contradicts transitivity of parallelism. Hence
$C\subseteq(i,j]\cup[j_0,i_0)$.
Then (2) gives the two asserted classes. Since there are no new
nonloops when $j<j_0$, the same assertion shows that every other
class remains a class.
\end{proof}

The preceding lemma determines the two resulting classes whenever a class splits. 
We then bound the number of classes that split.
Our approach involves the following technical lemma.
\begin{lem}\label{L: nested classes rectangle}
Suppose $a,b,c,d\in[n]$ are cyclically ordered, with $a\parallel d$, $b\parallel c$, and the two classes distinct. Let $b',c',d'$ be their lifts in $[a, a+n)$.
Then
\[
 \begin{split}
 &r[a,c]-r[b,c] + r[b,d]-r[a,d]\\
 &\qquad=|\{x \in [a,b'): \ c'<f(x)\le d'\}|\ge1.
 \end{split}
\]
\end{lem}
\begin{proof}
Apply~\eqref{EQ: interval rank 2} to the four integer intervals $[a,c']$, $[b',d']$, $[a,d']$, and $[b',c']$. Grouping them into two differences gives 
\[
 \begin{split}
 r[a,c]-r[b,c]
 &=|\{x \in [a,b') : \ f(x) >  c'\}|,\\
 r[a,d]-r[b,d]
 &=|\{x \in [a,b'): \ f(x) >  d'\}|.
 \end{split}
\]
Subtracting the second equality from the first gives the stated count.

It remains to prove that the count is positive. 
Let $X$ be a $k \times n$
matrix whose row span represents 
a general point of $\Pi_{u,z}$.
Let
$V[s,t]=\operatorname{span}\{X_x:x\in[s,t]\}$ for a cyclic interval
$[s,t]$. Since $b\parallel c$, equation~\eqref{EQ: pair rank} gives
$r[b,c]+r[c,b]=k+1$. Also, $[b,c]\cup[c,b]=[n]$, so
$V[b,c]+V[c,b]$ is the span of all columns of $X$, which has dimension $k$. 
Therefore
\[
 \dim\bigl(V[b,c]\cap V[c,b]\bigr)
 =r[b,c]+r[c,b]-k=1.
\]
The nonzero vector $X_b$ belongs to this one-dimensional space
$V[b,c]\cap V[c,b]$, so this space is exactly the span of $X_b$.
On the other hand, $a\in[c,b]$ implies $X_a \in V[c,b]$
while $X_a\notin V[b,c]\cap V[c,b]$ since $a$ is not parallel to $b$. 
Hence, $X_a\notin V[b,c]$.
The space $V[b,c]+\operatorname{span}(X_a)$ has dimension $r[b,c] + 1$.

Since $a\parallel d$, we have
\begin{equation}\label{EQ: V[b,c]}
 V[b,c]+\operatorname{span}(X_a)
 \subseteq V[a,c]\cap V[b,d],
\end{equation}
so the space on the right-hand side has dimension
at least $r[b,c] + 1$.
Notice that $V[a,c]+V[b,d]=V[a,d]$. It follows that
\[
 \begin{aligned}
 r[a,c]+r[b,d]
 &= \dim(V[a,c]) + \dim(V[b,d])\\
&= \dim\bigl(V[a,c]+V[b,d]\bigr)
   +\dim\bigl(V[a,c]\cap V[b,d]\bigr)\\
   &= \dim\bigl(V[a,d]\bigr)
   +\dim\bigl(V[a,c]\cap V[b,d]\bigr)
\ge r[a,d]+r[b,c]+1,
 \end{aligned}
\]
which proves the inequality.
\end{proof}

This lemma allows us to 
conclude either one or zero
class split. 

\begin{lem}\label{L: cover at most one split}
\label{L: cover opens loop preserves classes}\label{L: empty cover rectangle}
\begin{enumerate}
\item If $j<j_0$, at most one class in $\PP$ meets both $(i,j]$ and $[j_0,i_0)$.
\item If $j=j_0$, no class in $\PP$ meets both $(i,j]$ and $[j_0,i_0)$.
\end{enumerate}
\end{lem}

\begin{proof}
We first prove the claim
\begin{equation}\label{EQ: empty cover rectangle}
 \{x\in(i,j):f(x)\in(j_0,i_0+n)\}=\varnothing.
\end{equation}
Suppose the set has such an $x$ and let
\[
 q=z^{-1}(x),\qquad q_i=z^{-1}(i),\qquad q_j=z^{-1}(j).
\]
Since $z\lessdot_k t_{i,j}z$ and $i<x<j$, Lemma~\ref{L:cover-test}
excludes $q_i<q<q_j$. Thus either $q<q_i$ or $q>q_j$.
If $q<q_i\le k$, then
\[
 u(q)=f(x)-n<i_0=u(q_i),\qquad z(q)=x>i=z(q_i).
\]
This contradicts Proposition~\ref{P:k-Bruhat-endpoints}.
If $k<q_j<q$, then
\[
 u(q)=f(x)>j_0=u(q_j),\qquad z(q)=x<j=z(q_j),
\]
which contradicts Proposition~\ref{P:k-Bruhat-endpoints} again.

For~(1), suppose two distinct classes meet both $(i,j]$ and $[j_0,i_0)$.
Choose $a,b\in (i,j]$ and $c,d\in [j_0,i_0)$ with $a\parallel d$ and
$b\parallel c$ in the two classes. Without loss of generality, say $a<b$.
In the cyclic order starting at $a$, both $c$ and $d$ follow $b$.
The order $a,b,d,c$ is excluded by
Lemma~\ref{L: noncrossing parallel classes}. Hence $a,b,c,d$ are
cyclically ordered. 
Let $c', d'$ be the lifts
of $c,d$ in $[a, a+n)$.
We have
\[
 i<a<b\le j<j_0\le c'<d'<i_0+n.
\]
Lemma~\ref{L: nested classes rectangle} gives an integer
$x\in[a,b)$ with $f(x)\in(c',d']$. Since
$[a,b)\subseteq(i,j)$ and $(c',d']\subseteq(j_0,i_0+n)$, this
contradicts our claim~\eqref{EQ: empty cover rectangle}.

For~(2), suppose a class contains $a\in (i,j]$ and $b\in [j_0,i_0)$.
The number $j=j_0$ is a loop by Lemma~\ref{L: loop change}, whereas
$a,b$ are nonloops. The lift $b'$ of $b$ in $[a, a+n)$ therefore satisfies
$i<a<j<b'<i_0+n$.

Let $X$ be a $k\times n$ matrix whose row span represents
a general point of $\Pi_{u,z}$.
Recall that $X_x=X_{\overline{x}}$ for integer $x$.
Let $V_1,V_2,V_3$ be the spans of the columns indexed by
$[a,j)$, $(j,b']$, and $[a,b']$, respectively.
Since $j$ is a loop, $X_j=0$, so $V_1+V_2=V_3$.
Moreover, $V_1\cap V_2$ contains the span of $X_a$,
since $a\in[a,j)$, $b'\in(j,b']$, and $a\parallel b$.
Thus
\[
 r[a,j-1]+r[j+1,b']-r[a,b']
 =\dim(V_1\cap V_2)\ge1.
\]

On the other hand, applying~\eqref{EQ: interval rank}
to the three terms on the left gives
\[
\begin{aligned}
&r[a,j-1]+r[j+1,b']-r[a,b']\\
&\qquad=
|\{x\in[a,j):f(x)\in[j,b']\}|\\
&\qquad=
|\{x\in[a,j):f(x)\in(j,b']\}|,
\end{aligned}
\]
where the second equality follows from $f(j)=j$
and the injectivity of $f$.
The set is empty by~\eqref{EQ: empty cover rectangle},
since $[a,j)\subseteq(i,j)$ and
$(j,b']\subseteq(j_0,i_0+n)$.
This contradiction proves~(2).
\end{proof}

\begin{proof}[Proof of Proposition~\ref{P: general cover parallel classes}]
Suppose $j<j_0$. Lemma~\ref{L: loop change} gives $\LL_{u,z}=\LL_{u,z'}$.
By Lemma~\ref{L: cover rank rectangle}(2)--(3), each class meeting
both $(i,j]$ and $[j_0,i_0)$ splits into its two intersections, and all other
classes remain unchanged. Lemma~\ref{L: cover at most one split}(1)
shows that at most one class splits. This proves the first assertion.

Suppose $j=j_0$. Lemma~\ref{L: loop change} gives
$\LL_{u,z'}=\LL_{u,z}\setminus\{j\}$.
By Lemmas~\ref{L: cover rank rectangle}(2) and
\ref{L: cover opens loop preserves classes}(2), every old class remains
a class if we ignore $j$. Thus $j$ either joins one existing class or forms a singleton class.

For a nonloop $a$ of $\Pi_{u,z}$, equation~\eqref{EQ: class update rank} gives
\[
 r'[a,j]+r'[j,a]-r[a,j]-r[j,a]
 =\mathbf1_{\{a\notin[i_0,i]\}}.
\]
Since $j \in \LL_{u,z}$ and $a \notin \LL_{u,z}$,
\eqref{EQ: pair rank} gives $r[a,j]+r[j,a]-k\ge1$.
In $\Pi_{u,z'}$, both are nonloops. Hence they are parallel exactly when
\[
 i_0\le a\le i,\qquad r[a,j]+r[j,a]=k+1,
\]
that is, exactly when $a\in Q$. 
This proves the second assertion.
\end{proof}
\subsection{Proof of Proposition~\ref{P: second factor class update}}\label{SS: pair class update}
Return to the increasing chain~\eqref{EQ: chain from v to w}
and the integers $i'_s,j'_s$ defined at the beginning of the section.
We have $i'_s\le i_s<j_s\le j'_s$.
The chain is right-distinct, so position $z_{s-1}^{-1}(j_s)>k$
is not involved in an earlier swap of this chain. Hence
$v(z_{s-1}^{-1}(j_s))=j_s$.
The next lemma determines whether
$u$ and $v$
agree at this position.

\begin{lem}\label{L: j_s V(A)}
For $1\le s\le p$, exactly one of the following holds:
\begin{enumerate}
\item $j_s\in V(A)$, in which case $j_s<j'_s$,
$(j_s,j'_s)\in E(A)$, and $m(u,z_s)=m(u,z_{s-1})$;
\item $j_s\notin V(A)$, in which case $j'_s=j_s$ and
$m(u,z_s)=m(u,z_{s-1})+1$.
\end{enumerate}
\end{lem}
\begin{proof}
Put $q=z_{s-1}^{-1}(j_s)>k$. We have
$v(q)=j_s$ and $u(q)=j'_s$. If $u(q)>v(q)$, then
Lemma~\ref{L: first graph position values} gives
$(v(q),u(q))=(j_s,j'_s)\in E(A)$, so $j_s\in V(A)$.
If $u(q)=v(q)=j_s$, then Lemma~\ref{L: first graph position values} gives
$j_s\notin V(A)$.

Among positions greater than $k$, the permutations $z_{s-1},z_s$
differ only at $q$, where
\[
 z_{s-1}(q)=j_s,\qquad z_s(q)=i_s.
\]
If $j'_s>j_s$, then 
$u(q)>z_{s-1}(q)>z_s(q)$ shows that $q$ is counted by both
$m(u,z_{s-1})$ and $m(u,z_s)$.
If $j'_s=j_s$, it is counted only by $m(u,z_s)$.
Thus $m(u,z_s)=m(u,z_{s-1})$ in the first case and
$m(u,z_s)=m(u,z_{s-1})+1$ in the second.
\end{proof}

When $j_s\in V(A)$, Proposition~\ref{P: second factor class update}
asserts that $j_s,j'_s$ are parallel before the $s$th transposition.
We establish this first.

\begin{lem}\label{L: pair labels remain parallel}
If $j_s<j'_s$, then $j_s,j'_s$ belong to the same class of
$\PP_{u,z_{s-1}}$.
\end{lem}
\begin{proof}
By Lemmas~\ref{L: j_s V(A)} and~\ref{L: one factor classes},
$j_s,j'_s$ are parallel in $\Pi_{u,v}$.
Suppose they cease to be parallel before the $s$th transposition.
There is an integer $l<s$ such that they are in
one class of $\PP_{u,z_{l-1}}$ and in distinct classes of
$\PP_{u,z_l}$.
Proposition~\ref{P: general cover parallel classes}
forces exactly one of $j_s,j'_s$ to lie in $(i_l,j_l]$. 
Since $i_l<i_s<j_s<j'_s$, we have $j_s\le j_l<j'_s$.
Notice that $(i_l,j_l)$ and $(i_s,j_s)$ are two
edges of $B$, so $j_l\ne j_s$.
Combining these inequalities gives
$i_l<j_s<j_l<j'_s$.
This contradicts Lemma~\ref{L: mixed crossing reconstruction} applied
to $(j_s,j'_s)\in E(A)$ and $(i_l,j_l)\in E(B)$.
\end{proof}

\begin{proof}[Proof of Proposition~\ref{P: second factor class update}]
We argue by induction on $s$.
At $s=0$, Lemma~\ref{L: one factor classes} gives $|\PP_{u,v}|=k$,
and Proposition~\ref{P: k-Bruhat easy check} gives
$\ell(v)-\ell(u)=m(u,v)$. This proves~\eqref{EQ: two factor class equality} since $z_0=v$.

Now let $1\le s\le p$.
The inductive hypothesis is
\begin{equation}\label{EQ: two factor class equality 2}
 \ell(z_{s-1})-\ell(u)=m(u,z_{s-1})+|\PP_{u,z_{s-1}}|-k.
\end{equation}
We compare~\eqref{EQ: two factor class equality} with~\eqref{EQ: two factor class equality 2}.
Since $z_{s-1}\lessdot_k z_s$, the left-hand side increases by $1$.
It therefore suffices to prove
\begin{equation}\label{EQ: two factor class equality 3}
m(u,z_{s})+|\PP_{u,z_{s}}|=m(u,z_{s-1})+|\PP_{u,z_{s-1}}| + 1.
\end{equation}

We now consider the two cases in Lemma~\ref{L: j_s V(A)}. 
Suppose first that $j_s\in V(A)$, so $j_s < j'_s$
and $(j_s,j'_s)\in E(A)$.
Lemma~\ref{L: pair labels remain parallel} 
says $j_s,j'_s$ are in the same class $C \in \PP_{u,z_{s-1}}$.
They belong to $(i_s,j_s]$ and $[j'_s,i'_s)$, respectively.
Proposition~\ref{P: general cover parallel classes} says 
$\PP_{u,z_s}$ is obtained from $\PP_{u,z_{s-1}}$
by splitting $C$ into two classes,
so $|\PP_{u,z_{s}}| = |\PP_{u,z_{s-1}}| + 1$.
In addition, one new class is $C\cap (i_s,j_s]$,
so the other class is $C\setminus(i_s,j_s]$.
Lemma~\ref{L: j_s V(A)} gives $m(u,z_s)=m(u,z_{s-1})$,
so~\eqref{EQ: two factor class equality 3} follows.

Next, suppose $j_s\notin V(A)$.
Lemma~\ref{L: j_s V(A)} gives $j_s=j'_s$ and
$m(u,z_s)=m(u,z_{s-1})+1$.
Proposition~\ref{P: general cover parallel classes}
says $\PP_{u,z_s}$ is obtained from $\PP_{u,z_{s-1}}$
by adding $j_s$ to an existing class, or adding the singleton class $\{j_s\}$.
In the latter case, 
$|\PP_{u,z_{s}}| = |\PP_{u,z_{s-1}}| + 1$.
Combined with~\eqref{EQ: two factor class equality 2},
$m(u,z_s)+|\PP_{u,z_s}|-k$ would exceed
$\ell(z_s)-\ell(u)$, 
contradicting Proposition~\ref{P: numerical class bound}.
Thus $j_s$ joins an existing class as its largest member by Proposition~\ref{P: general cover parallel classes}. The number of classes is unchanged, 
so~\eqref{EQ: two factor class equality 3} follows.
\end{proof}

\section{The image of \texorpdfstring{$I$}{I}}\label{S: image of I}
Fix $u\le_k w$ with $\C(u,w,2)\ne\varnothing$, and put
$\rho=wu^{-1}$. We prove that every $I(u,v,w)$ is an ideal of
$\Gamma_{u,w}$ by describing its intersection with each $\Gamma_C$.
Recall that~\S\ref{sec:introduction} defines
\[
 \Gamma_{u,w}=\bigsqcup_{C\in\PP_{u,w}}\Gamma_C,\qquad
 \Gamma_C:\quad\rho^{-1}(c_0)\prec\cdots\prec\rho^{-1}(c_{L-1})
 \quad\text{where }C=\{c_0<\cdots<c_L\}.
\]
The last assertion of Lemma~\ref{L: finite parallel test} gives
$c_j<\rho^{-1}(c_j)\le c_{j+1}<\rho^{-1}(c_{j+1})$ for
$0\le j<L-1$. Thus the order on $\Gamma_C$ agrees with the usual
order on the integers.
We also use $\Gamma_C$ and $\Gamma_{u,w}$ for their underlying sets:
\begin{equation}
\label{EQ: Gamma as sets}
\begin{aligned}
\Gamma_C=\{\rho^{-1}(c):c\in C,\ c \ne \max C\},
\quad
\Gamma_{u,w}=\bigsqcup_{C\in\PP_{u,w}}\Gamma_C.
\end{aligned}
\end{equation}

An \definition{initial segment} of a finite set $M\subseteq[n]$
consists of its $t$ smallest members, for some $0\le t\le |M|$.
Clearly, if $M'\subseteq M$
\begin{equation}
\label{EQ: Initial segment intersection}
S \text{ is an initial segment of $M$} \Rightarrow     S\cap M' \text{ is an initial segment of $M'$}
\end{equation}

A subset $S\subseteq\Gamma_{u,w}$ is an ideal 
exactly when $S\cap\Gamma_C$ is an initial segment of $\Gamma_C$
for each $C\in\PP_{u,w}$.

The following proposition is the main result of this section.
It asserts that three sets are 
initial segments.

\begin{prop}\label{P: image ideal opening}
\label{P: class prefix and state}
Let $(u,v,w)\in\C(u,w,2)$, let $A$ be the increasing chain graph
from $u$ to $v$, and put $S=I(u,v,w)$.
For every $C=\{c_0<\cdots<c_L\}\in\PP_{u,w}$, there is an integer
$0\le t\le L$ such that
\begin{equation}\label{EQ: class prefix and state}
 \begin{aligned}
 V(A)\cap C&=\{c_0,\ldots,c_t\},\\
 \rho(S)\cap C&=\{c_0,\ldots,c_{t-1}\},\\
 S\cap\Gamma_C&=\{\rho^{-1}(c_0),\ldots,\rho^{-1}(c_{t-1})\}.
 \end{aligned}
\end{equation}
The last two sets are empty when $t=0$. Moreover,
$S\subseteq\Gamma_{u,w}$, so $S$ is an ideal of $\Gamma_{u,w}$.
\end{prop}

\begin{exa}\label{Ex: reconstruction class prefixes}
Continue Example~\ref{Ex: reconstruction scan},
\[
 u=[2,1,3,5,4],\qquad v=[2,5,1,3,4],\qquad w=[4,5,1,2,3],\qquad k=2,
\]
we have $S=I(u,v,w)=\{3\}$, $\rho(S)=\{1\}$
and $V(A) = \{1,2,3,5\}$.
Then we compute
\[
\begin{array}{c|c|c|c|c}
 C&\Gamma_C&t&V(A)\cap C&\rho(S)\cap C\\\hline
 \{1,5\}&\{3\}&1&\{1,5\}&\{1\}\\
 \{2\}&\varnothing&0&\{2\}&\varnothing\\
 \{3,4\}&\{4\}&0&\{3\}&\varnothing
\end{array}.
\]
For each $C$, we see that $S \cap \Gamma_C$
is an initial segment of $\Gamma_C$.
If we let $t$ be the size of $S \cap \Gamma_C$,
then $V(A) \cap C$ and $\rho(S) \cap C$
are initial segments of $C$ of sizes
$t+1$ and $t$ respectively. 
\end{exa}

Before giving the proof, we 
record the following consequences that will be used later.
\begin{lem}\label{L: Loops cap C}
Keep the notation of Proposition~\ref{P: image ideal opening}. Then
\[
 \LL_{u,v}\cap C=\{c_{t+1},\ldots,c_L\},\qquad
 \LL_{v,w}\cap C=\{c_0,\ldots,c_{t-1}\}.
\]
\end{lem}
\begin{proof}
By Lemma~\ref{L: first graph position values},
$[n] \setminus V(A)$ consists of $u(j)$
with $j > k$ and $u(j) = v(j)$.
Thus, this set is $\LL_{u,v}$.
By Proposition~\ref{P: image ideal opening},
\[
\LL_{u,v}\cap C = \big([n] \setminus V(A)\big) \cap C
= C \setminus \big(V(A)\cap C\big)=\{c_{t+1},\ldots,c_L\}.
\]

Next, we claim $\LL_{v,w}=\LL_{u,w}\sqcup\rho(S)$.
For $j > k$ with $v(j) = w(j)$, we have
either $u(j) = v(j)= w(j)$ or $u(j) > v(j)= w(j)$.
The former case is equivalent to $w(j) \in \LL_{u,w}$
while the latter case is equivalent to $w(j) \in \rho(S)$, proving our claim.
Using this claim and Proposition~\ref{P: image ideal opening},
\[
\LL_{v,w}\cap C = \big(\LL_{u,w}\sqcup\rho(S)\big) \cap C
= \rho(S) \cap C =\{c_0,\ldots,c_{t-1}\}. \qedhere
\]
\end{proof}

\begin{cor}\label{C: u v w length}
For $(u,v,w)\in\C(u,w,2)$, we have
\begin{align*}
 \LL_{u,v}&=\LL_{u,w}\sqcup
 \bigsqcup_{C\in\PP_{u,w}}\bigl(C\cap \LL_{u,v}\bigr),\\
 \LL_{v,w}&=\LL_{u,w}\sqcup
 \bigsqcup_{C\in\PP_{u,w}}\bigl(C\cap \LL_{v,w}\bigr).
\end{align*}
Consequently, if $S=I(u,v,w)$, then
\[
 \ell(v)-\ell(u)=|\PP_{u,w}|-k+|S|,\qquad
 \ell(w)-\ell(v)=m(u,w)-|S|.
\]
\end{cor}
\begin{proof}
If $j>k$ and $u(j)=w(j)$, then the inequalities
$u(j)\ge v(j)\ge w(j)$ give $u(j)=v(j)=w(j)$.
Hence $\LL_{u,w}$ is contained in both $\LL_{u,v}$ and $\LL_{v,w}$.
The parallel classes partition $[n]\setminus \LL_{u,w}$, so intersecting
these classes with $\LL_{u,v}$ or $\LL_{v,w}$ proves the
first two equalities.

Put $t_C=|S\cap\Gamma_C|$. Proposition~\ref{P: image ideal opening}
gives $|S|=\sum_Ct_C$. By Lemma~\ref{L: Loops cap C},
\begin{align*}
 |\LL_{u,v}|&=|\LL_{u,w}|+\sum_C\bigl(|C|-1-t_C\bigr)
          =n-|\PP_{u,w}|-|S|,\\
 |\LL_{v,w}|&=|\LL_{u,w}|+\sum_Ct_C=|\LL_{u,w}|+|S|.
\end{align*}
Then the first length identity in the corollary
follows from 
$\ell(v)-\ell(u)=n-k-|\LL_{u,v}|$.
The second follows from
$\ell(w)-\ell(v)=n-k-|\LL_{v,w}|$
and $m(u,w) = n - k - |\LL_{u,w}|$.
\end{proof}
\subsection{Proof of Proposition~\ref{P: image ideal opening}}
We start by describing $\Gamma_{u,w}$ as a set.
Put
\[
 R=\{r\in[n]:\rho(r)<r\}
   =\{u(j):j>k,\ u(j)>w(j)\}.
\]
\begin{lem}\label{L: Gamma membership}
We have $\Gamma_{u,w}\subseteq R$.
For $r\in R$, we have $\rho(r) \notin \LL_{u,w}$, and
$r\in\Gamma_{u,w}$ if and only if $\rho(r)$ is not the maximum of
its parallel class.
\end{lem}
\begin{proof}
Every member of $\Gamma_{u,w}$ is $\rho^{-1}(a)$ for a nonmaximal
member $a$ of a parallel class. Choose $b>a$ in that class. The last
assertion of Lemma~\ref{L: finite parallel test} gives
$a<\rho^{-1}(a)\le b$, so $\rho^{-1}(a)\in R$.

Now take $r\in R$. If $\rho(r)$ is a loop, then
Lemma~\ref{L: loops endpoints} says
$\rho(\rho(r))=\rho(r)$, so $r=\rho(r)$, contrary to $r\in R$. 
Thus $\rho(r)$ belongs to a parallel class. 
The last statement in the lemma follows from~\eqref{EQ: Gamma as sets}.
\end{proof}

Fix $(u,v,w)\in\C(u,w,2)$, with increasing chain graphs $A,B$,
and write its second chain as in~\eqref{EQ: chain from v to w}.
Put
\[
 U_s=\{j_1,\ldots,j_s\}\quad(1\le s\le p),\qquad U_0=\varnothing,
 \qquad S=I(u,v,w).
\]
By Lemma~\ref{L: first graph position values},
$E(B) = \{(i_s, j_s): s \in [p] \}$.
Thus, the same lemma implies
\begin{equation}
\label{EQ: U_p}
U_p=\{j_1,\ldots,j_p\}
    =\{b:(a,b)\in E(B)\}
    =\{v(j):j>k,\ w(j)<v(j)\}.
\end{equation}

The preceding section determines how the parallel classes change
in each step along the increasing chain from $v$ to $w$.
The next lemma uses this description.

\begin{lem}\label{L: class marked maximum}
For $0\le s\le p$ and $C\in\PP_{u,z_s}$, the following hold.
\begin{enumerate}
\item The set $C\cap V(A)$ is a nonempty initial segment of $C$.
\item Its only member in $\Tail(A)\cup U_s$ is $\max(C\cap V(A))$.
\end{enumerate}

\end{lem}
\begin{proof}
We argue by induction on $s$.
At $s=0$, Lemma~\ref{L: one factor classes} identifies each class
of $\PP_{u,v}$ with an ordered path in $A$.
Thus $C\subseteq V(A)$, proving~(1). The ordered path 
has exactly one element from $\Tail(A)$: its largest number.  
This proves~(2).

Suppose $s>0$ and the lemma holds at $s-1$.
The set $U_s$ differs from $U_{s-1}$ by adding $j_s$.
If $C$ is also a class in $\PP_{u, z_{s-1}}$ and $j_s\notin C$, 
the inductive step follows from the induction hypothesis. 
It remains to assume $C \notin \PP_{u, z_{s-1}}$ 
or $j_s \in C$.

First suppose $j_s\notin V(A)$.
By Proposition~\ref{P: second factor class update}, the only 
class not in $\PP_{u,z_{s-1}}$ 
has the form $C=C'\cup\{j_s\}$, where
$C'\in\PP_{u,z_{s-1}}$ and $j_s = \max(C)$.
We have $V(A)\cap C=V(A)\cap C'$ which is still a nonempty
initial segment of $C$, proving~(1).
Also, $U_s\cap V(A)=U_{s-1}\cap V(A)$, proving~(2).

Now suppose $j_s\in V(A)$. 
Proposition~\ref{P: second factor class update}
says the class $D \in \PP_{u, z_{s-1}}$
containing $j_s,j'_s$ splits into two classes in $\PP_{u, z_s}$.
They are $D_1 = D\cap(i_s,j_s]$ and $D_2 = D\setminus(i_s,j_s]$.
The edge $(j_s,j'_s)\in E(A)$ shows that $j_s\in D_1\cap V(A)$
and $j'_s\in D_2\cap V(A)$.
Thus, both new classes are nonempty. By induction, $D\cap V(A)$
is an initial segment of $D$. 
By~\eqref{EQ: Initial segment intersection}, 
$D_1 \cap V(A)$ and $D_2 \cap V(A)$
are nonempty initial segments of $D_1$ and $D_2$ respectively.

Put $d=\max(D\cap V(A))$. By the inductive hypothesis,
\begin{equation}
\label{EQ: D V(A)}  
 D\cap V(A)\cap\bigl(\Tail(A)\cup U_{s-1}\bigr)=\{d\}.
\end{equation}
Since $j_s<j'_s\le d$, we have $d\in D_2$, and
\[
 \max(D_1\cap V(A))=j_s,\qquad \max(D_2\cap V(A))=d.
\]
From~\eqref{EQ: D V(A)},
$D_1\cap V(A)$ contains no member in
$\Tail(A)\cup U_{s-1}$. 
Thus, its unique member of $\Tail(A) \cup U_s$
is $j_s$.
For $D_2\cap V(A)$, its unique member in
$\Tail(A)\cup U_{s-1}$ is $d$.
Since $j_s \notin D_2$,
the unique member of $D_2\cap V(A)$ in 
$\Tail(A)\cup U_{s}$ is still $d$.
\end{proof}

The lemma above studies $V(A) \cap C$.
As outlined by Proposition~\ref{P: image ideal opening}, 
we study $\rho(S)$ next. 

\begin{lem}\label{L: state settled vertices}
We have 
\begin{equation}\label{EQ: state settled vertices}
 \rho(S)=V(A)\setminus\bigl(\Tail(A)\cup U_p\bigr).
\end{equation}
\end{lem}
\begin{proof}
Recall from~\eqref{EQ: U_p}
that $U_p$ consists of $v(j)$ with $j>k$ and $w(j)<v(j)$.
Take $a\in V(A)\setminus(\Tail(A)\cup U_p)$.
Since $a\notin\Tail(A)$, there is a unique edge $(a,r)\in E(A)$.
Let $q = u^{-1}(r)$.
Lemma~\ref{L: first graph position values} gives
$q>k$ and $v(q)=a<r$.
Since $a\notin U_p$, $w(q)=v(q) = a$.
Thus, $\rho(r)=a$ and $r\in S$.

Conversely, take $r\in S$ and put $q=u^{-1}(r)$ and $a=\rho(r)$.
By the definition of $S$, $v(q)=w(q)=a<r=u(q)$.
Proposition~\ref{P:k-Bruhat-endpoints} gives $q>k$.
Lemma~\ref{L: first graph position values} gives
$(a,r)\in E(A)$. Thus $a\in V(A) \setminus \Tail(A)$.
Since $v(q)=w(q) = a$, $a \notin U_p$,
proving $\rho(r) = a \in V(A)\setminus\bigl(\Tail(A)\cup U_p\bigr)$.
\end{proof}

\begin{proof}[Proof of Proposition~\ref{P: image ideal opening}]
Fix $C=\{c_0<\cdots<c_L\}\in\PP_{u,w}$.
Lemma~\ref{L: class marked maximum} at $z_p=w$ gives an integer $0\le t\le L$ such that
\[
 C\cap V(A)=\{c_0,\ldots,c_t\},\qquad
 C\cap V(A)\cap(\Tail(A)\cup U_p)=\{c_t\}.
\]
Lemma~\ref{L: state settled vertices} then gives
\[
\rho(S)\cap C
=(V(A)\cap C)\setminus\bigl(\Tail(A)\cup U_p\bigr).
\]
The only element in $V(A) \cap C$ that lies in $\Tail(A)\cup U_p$
is $c_t$, so the right hand side becomes $\{c_0,\ldots,c_{t-1}\}$.
Notice that $\rho(S)$ cannot contain the largest number of the class $C$.

To show $S \subseteq \Gamma_{u,w}$,
take $r\in S$. By definition, $\rho(r)<r$, so $r\in R$.
Lemma~\ref{L: Gamma membership} shows that $\rho(r)$ is a nonloop.
It therefore belongs to some class in $\PP_{u,w}$.
The preceding paragraph shows that it is not the maximum of
its class, so Lemma~\ref{L: Gamma membership} 
gives $r \in \Gamma_{u,w}$.
Thus $S\subseteq\Gamma_{u,w}$.

Finally, $\rho(\Gamma_C)=C\setminus\{c_L\}$, and $c_L\notin\rho(S)\cap C$.
Taking $\rho^{-1}$ in the second equality gives the third equality
in~\eqref{EQ: class prefix and state}. Thus $S$ meets every
$\Gamma_C$ in an initial segment and is an ideal of $\Gamma_{u,w}$.
\end{proof}

\section{Surjectivity of \texorpdfstring{$I$}{I}}\label{S: surjectivity}
Fix $u\le_k w$ with $\C(u,w,2)\ne\varnothing$, and put $\rho=wu^{-1}$.
Write $\J$ for the set of ideals of $\Gamma_{u,w}$.
Proposition~\ref{P: image ideal opening} gives $\img(I)\subseteq\J$,
and Proposition~\ref{P: reconstruction} shows that $I$ is injective.
It remains to prove that every ideal occurs in the image.

\begin{prop}\label{P: surjectivity}
The map $I:\C(u,w,2)\longrightarrow\J$ is surjective. 
\end{prop}

We introduce a useful convention that will be used
throughout this section. 
For a nonempty finite set $C$ of integers and $0\le t<|C|$, let $C_t$
be its $t$ smallest members and let $C^t$ be its $|C|-t-1$ largest
members. Thus, if $c$ is the $(t+1)$st smallest member of $C$, then
$C_t=\{a\in C:a<c\}$ and $C^t=\{a\in C:a>c\}$.

\begin{rem}\label{R: identify ideal with maps}
An ideal $S\in\J$ is uniquely determined by the integers
$|S\cap\Gamma_C|$ for each $C \in \PP_{u,w}$.
We identify $S$ with the map $\bt$
sending each $C \in \PP_{u,w}$
to $|S\cap\Gamma_C|$.
In this way, $\J$ is identified with the set of maps
defined on $\PP_{u,w}$,
sending $C$ to a number in $\{0, 1, \dots, |C| - 1\}$.
Conversely, the ideal associated to such a map $\bt$ 
is
\[
 S_{\bt} :=\bigsqcup_{C\in\PP_{u,w}}\rho^{-1}(C_{\bt(C)}),
\]
which implies $\rho(S_{\bt})=\bigsqcup_C C_{\bt(C)}$.
\end{rem}

\begin{exa}
Suppose $\PP_{u,w} = \{C, D\}$
where $\Gamma_C$ is $c_1\prec c_2$ and $\Gamma_D$ 
has only $d$.
All ideals of $\Gamma_{u,w}$ are listed below. 
We write down their corresponding $\bt$ underneath
as a pair $(\bt(C), \bt(D))$.
\[
\begin{array}{c|cccccc}
 S&\varnothing&\{c_1\}&\{c_1,c_2\}&\{d\}&\{c_1,d\}&\{c_1,c_2,d\}\\
 (\bt(C),\bt(D))&(0,0)&(1,0)&(2,0)&(0,1)&(1,1)&(2,1).
\end{array} \qedhere
\]
\end{exa}

Under this identification, we also write $\bt\in\J$ and
$\bt\in\img(I)$. The indicator $\mathbf1_{\{\bt\in\img(I)\}}$ is
$1$ when the ideal associated to $\bt$ belongs to the image of $I$, and $0$ otherwise.
Since $I$ is injective, this counts $(u,v,w) \in \C(u,w,2)$ mapped to $\bt$.

Recall from Subsection~\ref{SS:double-weight} that $\y=(y_1,\ldots,y_n)$ denotes the $y$-variables.
For $C=\{c_0<\cdots<c_{L_C}\}$ and $0\le t\le L_C$, define
\[
 \phi_{C,t}(x_1,x_2;\y)
 =\prod_{j \in C_t}(x_2-y_{j})
  \prod_{j \in C^t}(x_1-y_j).
\]
We use Corollary~\ref{C: u v w length} to rewrite the polynomial
$F'_{u,w}(x_1,x_2;\y)$ from~\eqref{EQ: double chain weight} as a sum over $\J$.

\begin{lem}
We have
\begin{equation}\label{EQ: occurrence weighted sum}
 F'_{u,w}(x_1,x_2;\y)
 =H(x_1,x_2;\y)\prod_{a\in \LL_{u,w}}(x_1-y_a)(x_2-y_a),
\end{equation}
\[
 \text{ where } H(x_1,x_2;\y)=\sum_{\bt\in\J}
 \mathbf1_{\{\bt\in\img(I)\}}
 \prod_{C\in\PP_{u,w}}\phi_{C,\bt(C)}(x_1,x_2;\y).
\]
\end{lem}
\begin{proof}
By definition, the contribution of $(u,v,w)\in\C(u,w,2)$ to $F'_{u,w}$ is
\[
 \prod_{j\in \LL_{u,v}}(x_1-y_j)\prod_{j\in \LL_{v,w}}(x_2-y_j).
\]
The two set equalities in Corollary~\ref{C: u v w length} rewrite it as
\[
 \prod_{j\in \LL_{u,w}}(x_1-y_j)(x_2-y_j)
 \prod_{C\in\PP_{u,w}}
 \left(\prod_{j\in C\cap \LL_{u,v}}(x_1-y_j)
       \prod_{j\in C\cap \LL_{v,w}}(x_2-y_j)\right).
\]
Identify $I(u,v,w)$ with $\bt$ as in
Remark~\ref{R: identify ideal with maps}.
Write $C = \{ c_0 < c_1 < \cdots < c_L\}$.
Lemma~\ref{L: Loops cap C} gives
\[
 C\cap \LL_{u,v}=\{c_{\bt(C)+1},\ldots,c_{L}\} = C^{\bt(C)},\qquad
 C\cap \LL_{v,w}=\{c_0,\ldots,c_{\bt(C)-1}\} = C_{\bt(C)}.
\]
Thus the contribution from $C$ is exactly 
$\phi_{C,\bt(C)}(x_1,x_2;\y)$.
Summing over $(u,v,w)$ and grouping the contributions by $\bt\in\J$
proves~\eqref{EQ: occurrence weighted sum}.
\end{proof}

The symmetry of $F'_{u,w}$ gives the symmetry needed below.
\begin{cor}\label{C: P sym}
The polynomial $H(x_1,x_2;\y)$ is symmetric in $x_1,x_2$.
\end{cor}
\begin{proof}
By Proposition~\ref{P: double weight symmetry},
$F'_{u,w}(x_1,x_2;\y)$ is symmetric in $x_1,x_2$.
The product over $\LL_{u,w}$ in~\eqref{EQ: occurrence weighted sum} is
nonzero and symmetric. Canceling it proves the claim.
\end{proof}

We next show that changing the size of an ideal in one chain by one
cannot change whether it belongs to $\img(I)$.
\begin{lem}\label{L: neighbouring ideal coefficients}
Suppose $\bt_1,\bt_2\in\J$ differ by $1$ on one class $C\in\PP_{u,w}$
and agree on every other class. Then
\[
 \bt_1\in\img(I)\quad\Longleftrightarrow\quad\bt_2\in\img(I).
\]
\end{lem}
\begin{proof}
Write $C=\{c_0<\cdots<c_L\}$. Without loss of generality,
assume $\bt_1(C)=i$ and $\bt_2(C)=i+1$, where $0\le i<L$.

For every $D\ne C$, write $D=\{d_0<\cdots<d_{L_D}\}$.
Consider the set 
$$
(C \setminus \{c_{i}, c_{i+1}\}) \sqcup \bigsqcup_{\substack{D\in\PP_{u,w}\\D\ne C}} \big( D_{\bt_1(D)} \sqcup D^{\bt_1(D)} \big).
$$

It corresponds to a squarefree monomial $Y$ of $y$-variables. 
We extract the coefficient of $Y$ from $H$ 
viewed as a polynomial in $\y$ with coefficients in $\mathbb{Z}[x_1,x_2]$.
For $D\in \PP_{u,w}$, 
let $D(t)$ denote the $t^\text{th}$ 
smallest number of $D$.
For $D \neq C$,
$Y$ contains every $y_j$ for
$j\in D$ except $y_{D(\bt_1(D) + 1)}$.
On the other hand, $\phi_{D,t}$ involves 
all variables except $y_{D(t+1)}$.
Thus, a summand $\bt$ can contribute only when $\bt(D)=\bt_1(D)$.
Similarly, the only possibilities for $\bt(C)$
are $i$ and $i+1$. Hence
\[
 [Y]H(x_1, x_2; \y)=(-1)^{\deg Y}
 \left(\mathbf1_{\{\bt_1\in\img(I)\}}x_1
       +\mathbf1_{\{\bt_2\in\img(I)\}}x_2\right).
\]
By Corollary~\ref{C: P sym}, this coefficient is symmetric in $x_1,x_2$.
The two indicator functions are therefore equal, proving the assertion.
\end{proof}

\begin{proof}[Proof of Proposition~\ref{P: surjectivity}]
View $\J$ as maps on $\PP_{u,w}$ as in
Remark~\ref{R: identify ideal with maps}. 
These maps are connected by the operator
that changes the value on one class by $1$.
Lemma~\ref{L: neighbouring ideal coefficients} shows that each
application of this operator preserves membership in $\img(I)$. Hence $\img(I)$ is
empty or equals $\J$.
Since $\C(u,w,2)\ne\varnothing$, Proposition~\ref{P: image ideal opening}
excludes the empty case. Thus $\img(I)=\J$.
\end{proof}

\section{Proofs of the main results}\label{S: main proofs}
We now prove the theorems stated in the introduction by combining
the results of the preceding sections. 

\begin{proof}[Proof of Theorem~\ref{T: main}]
Since $m(u,w)=n-k-|\LL_{u,w}|$, the inequality follows from
Proposition~\ref{P: numerical class bound} and
Corollary~\ref{C: two factor equality}.
Assume equality holds.

 Propositions~\ref{P: image ideal opening},
\ref{P: reconstruction}, and~\ref{P: surjectivity} show that $I$ is a
bijection to the ideals of $\Gamma_{u,w}$ with inverse $\Chain(\cdot)$.

The parallel classes partition $[n]\setminus\LL_{u,w}$, 
so the sum of entries in $\alpha$ gives $m(u,w)-d$. Corollary~\ref{C: u v w length} 
therefore gives
\[
 \ell(v)-\ell(u)=d+|S|,\qquad
 \ell(w)-\ell(v)=d+|\alpha|-|S|.
\]
An ideal independently selects an initial segment of each $\Gamma_C$,

\[
 F_{u,w}(x_1,x_2)
 =(x_1x_2)^d\prod_{C\in\PP_{u,w}}
 \sum_{t=0}^{|C|-1}x_1^t x_2^{|C|-1-t}
 =(x_1x_2)^d h_\alpha(x_1,x_2).
 \qedhere
\]
\end{proof}

\begin{proof}[Proof of Theorem~\ref{T: two row coefficients}]
Follows from Theorem~\ref{T: main}.

\end{proof}

\begin{proof}[Proof of Theorem~\ref{T: double two row}]
If~\eqref{EQ: introduction numerical bound} is not 
equality, Theorem~\ref{T: main} gives
$\C(u,w,2)=\varnothing$, so $F'_{u,w}=0$.
Otherwise, Proposition~\ref{P: surjectivity} makes every indicator function in
\eqref{EQ: occurrence weighted sum} equal to $1$.
We can then factor $H(x_1, x_2; \y)$ 
over the parallel classes.
For $C=\{c_0<\cdots<c_L\}$, notice that 
\[
 \sum_{t=0}^{L}\phi_{C,t}(x_1,x_2;\y)
 =\frac{\prod_{c\in C}(x_1-y_c)-\prod_{c\in C}(x_2-y_c)}{x_1-x_2}.
\]
Substituting into~\eqref{EQ: occurrence weighted sum} gives the proof.

\end{proof}

\bibliographystyle{alpha}
\bibliography{pieri}

@article{LS:Schubert,
  author  = {Lascoux, Alain and Sch{\"u}tzenberger, Marcel-Paul},
  title   = {Polyn{\^o}mes de {Schubert}},
  journal = {C. R. Acad. Sci. Paris S{\'e}r. I Math.},
  volume  = {294},
  number  = {13},
  year    = {1982},
  pages   = {447--450},
  url     = {https://www-igm.univ-mlv.fr/~berstel/Mps/Travaux/A/1982-1SchubertCras.pdf}
}

@article{Monk,
  author  = {Monk, D.},
  title   = {The geometry of flag manifolds},
  journal = {Proc. London Math. Soc. (3)},
  volume  = {9},
  number  = {2},
  year    = {1959},
  pages   = {253--286},
  doi     = {10.1112/plms/s3-9.2.253}
}

@article{Sot,
  author  = {Sottile, Frank},
  title   = {{Pieri}'s formula for flag manifolds and {Schubert} polynomials},
  journal = {Ann. Inst. Fourier (Grenoble)},
  volume  = {46},
  number  = {1},
  year    = {1996},
  pages   = {89--110},
  doi     = {10.5802/aif.1508},
  url     = {https://numdam.org/articles/10.5802/aif.1508/}
}

@article{BS:monoid,
  author  = {Bergeron, Nantel and Sottile, Frank},
  title   = {A monoid for the {Grassmannian Bruhat} order},
  journal = {European J. Combin.},
  volume  = {20},
  number  = {3},
  year    = {1999},
  pages   = {197--211},
  doi     = {10.1006/eujc.1999.0283},
  eprint  = {math/9712258},
  archivePrefix = {arXiv},
  url     = {https://franksottile.github.io/research/pdf/monoid.pdf}
}

@article{BS:skewSchubert,
  author  = {Bergeron, Nantel and Sottile, Frank},
  title   = {Skew {Schubert} functions and the {Pieri} formula for flag manifolds},
  journal = {Trans. Amer. Math. Soc.},
  volume  = {354},
  number  = {2},
  year    = {2002},
  pages   = {651--673},
  doi     = {10.1090/S0002-9947-01-02845-8},
  eprint  = {alg-geom/9709034},
  archivePrefix = {arXiv},
  url     = {https://arxiv.org/abs/alg-geom/9709034}
}

@article{Hall:representatives,
  author  = {Hall, P.},
  title   = {On representatives of subsets},
  journal = {J. London Math. Soc.},
  volume  = {10},
  number  = {1},
  year    = {1935},
  pages   = {26--30},
  doi     = {10.1112/jlms/s1-10.37.26}
}

@article{KLS:positroids,
  author  = {Knutson, Allen and Lam, Thomas and Speyer, David E.},
  title   = {Positroid varieties: juggling and geometry},
  journal = {Compos. Math.},
  volume  = {149},
  number  = {10},
  year    = {2013},
  pages   = {1710--1752},
  doi     = {10.1112/S0010437X13007240},
  eprint  = {1111.3660},
  archivePrefix = {arXiv},
  primaryClass = {math.AG},
  url     = {https://arxiv.org/abs/1111.3660}
}

@misc{Postnikov:total,
  author  = {Postnikov, Alexander},
  title   = {Total positivity, {Grassmannians}, and networks},
  year    = {2006},
  howpublished = {Preprint, arXiv:math/0609764},
  eprint  = {math/0609764},
  archivePrefix = {arXiv},
  primaryClass = {math.CO},
  url     = {https://arxiv.org/abs/math/0609764}
}

@article{Oh:positroids,
  author  = {Oh, Suho},
  title   = {Positroids and {Schubert} matroids},
  journal = {J. Combin. Theory Ser. A},
  volume  = {118},
  number  = {8},
  year    = {2011},
  pages   = {2426--2435},
  doi     = {10.1016/j.jcta.2011.06.006},
  eprint  = {0803.1018},
  archivePrefix = {arXiv},
  primaryClass = {math.CO},
  url     = {https://arxiv.org/abs/0803.1018}
}

@article{MO:rank,
  author  = {Mcalmon, Robert and Oh, Suho},
  title   = {The rank function of a positroid and non-crossing partitions},
  journal = {Electron. J. Combin.},
  volume  = {27},
  number  = {1},
  year    = {2020},
  pages   = {P1.11},
  note    = {13 pp.},
  doi     = {10.37236/8256},
  eprint  = {1709.01580},
  archivePrefix = {arXiv},
  primaryClass = {math.CO},
  url     = {https://www.combinatorics.org/ojs/index.php/eljc/article/view/v27i1p11}
}

@article{Samuel:Pieri,
  author  = {Samuel, Matthew J.},
  title   = {A {Molev--Sagan} type formula for double {Schubert} polynomials},
  journal = {J. Pure Appl. Algebra},
  volume  = {228},
  number  = {7},
  year    = {2024},
  pages   = {107636},
  doi     = {10.1016/j.jpaa.2024.107636},
  eprint  = {2401.11060},
  archivePrefix = {arXiv},
  primaryClass = {math.CO},
  url     = {https://arxiv.org/abs/2401.11060}
}

@misc{CY:doubleEG,
  author  = {Chou, Jack Chen-An and Yu, Tianyi},
  title   = {A positive combinatorial formula for the double {Edelman--Greene} coefficients},
  year    = {2025},
  howpublished = {Preprint, arXiv:2512.12185},
  eprint  = {2512.12185},
  archivePrefix = {arXiv},
  primaryClass = {math.CO},
  url     = {https://arxiv.org/abs/2512.12185}
}

@article{Pawlowski:positroids,
  author  = {Pawlowski, Brendan},
  title   = {A representation-theoretic interpretation of positroid classes},
  journal = {Adv. Math.},
  volume  = {429},
  year    = {2023},
  pages   = {109178},
  doi     = {10.1016/j.aim.2023.109178},
  eprint  = {1612.00097},
  archivePrefix = {arXiv},
  primaryClass = {math.CO},
  url     = {https://arxiv.org/abs/1612.00097}
}

@article{Kogan:RC,
 author={Kogan, Mikhail},
 title={{RC}-graphs and a generalized {Littlewood--Richardson} rule},
 journal={Int. Math. Res. Not.}, year={2001}, number={15}, pages={765--782},
 doi={10.1155/S1073792801000393},
 eprint={math/0010108}, archivePrefix={arXiv}
}

@article{Lenart:growth,
 author={Lenart, Cristian},
 title={Growth diagrams for the {Schubert} multiplication},
 journal={J. Combin. Theory Ser. A}, year={2010}, volume={117}, number={7},
 pages={842--856}, doi={10.1016/j.jcta.2010.02.004},
 eprint={0901.4149}, archivePrefix={arXiv}
}

@article{Huang:separated,
 author={Huang, Daoji},
 title={{Schubert} products for permutations with separated descents},
 journal={Int. Math. Res. Not. IMRN}, year={2023}, number={20},
 pages={17461--17493}, doi={10.1093/imrn/rnac299},
 eprint={2105.01591}, archivePrefix={arXiv}
}

@misc{KZ:separated,
 author={Knutson, Allen and Zinn-Justin, Paul},
 title={{Schubert} puzzles and integrability {III}: separated descents},
 year={2023}, howpublished={Preprint, arXiv:2306.13855},
 eprint={2306.13855}, archivePrefix={arXiv}, primaryClass={math.CO},
 url={https://arxiv.org/abs/2306.13855}
}

@misc{TaoZheng:twoRow,
 author={Tao, Zijie and Zheng, Yunchi},
 title={Bounds on {Schubert} coefficients in the two-row case},
 year={2024}, howpublished={Preprint, arXiv:2411.19735},
 eprint={2411.19735}, archivePrefix={arXiv}, primaryClass={math.CO},
 url={https://arxiv.org/abs/2411.19735}
}

@incollection{FK:quadratic,
 author={Fomin, Sergey and Kirillov, Anatol N.},
 title={Quadratic algebras, {Dunkl} elements, and {Schubert} calculus},
 booktitle={Advances in Geometry},
 editor={Brylinski, Jean-Luc and Brylinski, Ranee and Nistor, Victor and Tsygan, Boris and Xu, Ping},
 series={Progress in Mathematics}, volume={172},
 publisher={Birkh{\"a}user}, address={Boston, MA}, year={1999},
 pages={147--182}, doi={10.1007/978-1-4612-1770-1_8}
}

@incollection{Postnikov:quantumPieri,
 author={Postnikov, Alexander},
 title={On a quantum version of {Pieri}'s formula},
 booktitle={Advances in Geometry},
 editor={Brylinski, Jean-Luc and Brylinski, Ranee and Nistor, Victor and Tsygan, Boris and Xu, Ping},
 series={Progress in Mathematics}, volume={172},
 publisher={Birkh{\"a}user}, address={Boston, MA}, year={1999},
 pages={371--383}, doi={10.1007/978-1-4612-1770-1_15},
 url={https://math.mit.edu/~apost/papers/pieri.pdf}
}

@article{MPP:SchurSchubert,
 author={M{\'e}sz{\'a}ros, Karola and Panova, Greta and Postnikov, Alexander},
 title={{Schur} times {Schubert} via the {Fomin--Kirillov} algebra},
 journal={Electron. J. Combin.}, volume={21}, number={1}, year={2014},
 pages={P1.39}, doi={10.37236/3659},
 eprint={1210.1295}, archivePrefix={arXiv}
}

@article{BS:BruhatGeometry,
 author={Bergeron, Nantel and Sottile, Frank},
 title={{Schubert} polynomials, the {Bruhat} order, and the geometry of flag manifolds},
 journal={Duke Math. J.}, volume={95}, number={2}, year={1998},
 pages={373--423}, doi={10.1215/S0012-7094-98-09511-4},
 eprint={alg-geom/9703001}, archivePrefix={arXiv}
}

@article{LenartSottile:skew,
 author={Lenart, Cristian and Sottile, Frank},
 title={Skew {Schubert} polynomials},
 journal={Proc. Amer. Math. Soc.}, volume={131}, number={11}, year={2003},
 pages={3319--3328}, doi={10.1090/S0002-9939-03-06919-3},
 eprint={math/0202090}, archivePrefix={arXiv}
}

@article{LenartSottile:KPieri,
 author={Lenart, Cristian and Sottile, Frank},
 title={A {Pieri}-type formula for the {$K$}-theory of a flag manifold},
 journal={Trans. Amer. Math. Soc.}, volume={359}, number={5}, year={2007},
 pages={2317--2342}, doi={10.1090/S0002-9947-06-04043-8},
 eprint={math/0407412}, archivePrefix={arXiv}
}

@article{ARW:noncrossing,
 author={Ardila, Federico and Rinc{\'o}n, Felipe and Williams, Lauren},
 title={Positroids and non-crossing partitions},
 journal={Trans. Amer. Math. Soc.}, volume={368}, number={1}, year={2016},
 pages={337--363}, doi={10.1090/tran/6331},
 eprint={1308.2698}, archivePrefix={arXiv}
}

@article{GalashinLam:cluster,
 author={Galashin, Pavel and Lam, Thomas},
 title={Positroid varieties and cluster algebras},
 journal={Ann. Sci. {\'E}c. Norm. Sup{\'e}r. (4)},
 volume={56}, number={3}, year={2023}, pages={859--884},
 doi={10.24033/asens.2545}, eprint={1906.03501}, archivePrefix={arXiv}
}
\end{document}